\documentclass{amsart}

\usepackage{amssymb}
\usepackage{amsmath}
\usepackage{amsthm}
\usepackage{amsfonts}

\usepackage{a4wide}
\usepackage{longtable}
\usepackage{multirow}

\newtheorem{lemma}{Lemma}

\newtheorem{theorem}{Theorem}
\newtheorem{corollary}{Corollary}

\theoremstyle{definition}
\newtheorem{definition}{Definition}
\newtheorem{example}{Example}

\numberwithin{equation}{section}

\begin{document}

%\linenumbers

\title[Schur \(\sigma\)-Groups of Scholz-Taussky Type F]
{Schur \(\sigma\)-Groups of Scholz-Taussky Type F}

\author{Daniel C. Mayer}
\address{Naglergasse 53\\8010 Graz\\Austria}
\email{algebraic.number.theory@algebra.at}
\urladdr{http://www.algebra.at}

\thanks{Research supported by the Austrian Science Fund (FWF): projects J0497-PHY, P26008-N25, and by EUREA}



\date{Sunday, 14 August 2022}

%--------------------------------------------------------------------------------

\begin{abstract}
For finite metabelian \(3\)-groups \(M\)
with elementary bicyclic commutator quotient \(M/M^\prime\simeq C_3\times C_3\),
coclass \(\mathrm{cc}(M)\in\lbrace 4,6\rbrace\),
and transfer kernel type \(\mathrm{F}\),
the smallest Schur \(\sigma\)-groups \(S\)
with second derived quotient \(S/S^{\prime\prime}\simeq M\)
are determined.
Evidence is provided of arithmetical realizations of these groups
by second \(3\)-class groups
\(M\simeq\mathrm{G}_3^2(K)=\mathrm{Gal}(\mathrm{F}_3^2(K)/K)\),
respectively \(3\)-class field tower groups
\(S\simeq\mathrm{G}_3^\infty(K)=\mathrm{Gal}(\mathrm{F}_3^\infty(K)/K)\),
of imaginary quadratic number fields \(K=\mathbb{Q}(\sqrt{d})\).
\end{abstract}

\maketitle

%\newpage
%--------------------------------------------------------------------------------

\section{Introduction}
\label{s:Intro}

\noindent
A \textit{Schur \(\sigma\)-group} \(S\)
has a balanced presentation \(d(S)=r(S)\) with coinciding
generator rank  \(d(S)=\dim_{\mathbb{F}_p}H^1(S,\mathbb{F}_p)\) and
relation rank \(r(S)=\dim_{\mathbb{F}_p}H^2(S,\mathbb{F}_p)\),
and an automorphism \(\sigma\in\mathrm{Aut}(S)\)
inducing the inversion \(x\mapsto x^{-1}\)
on the first and second cohomology group,
\(H^1(S,\mathbb{F}_p)\) and \(H^2(S,\mathbb{F}_p)\).
Therefore, \(\sigma\) is called \textit{generator- and relator-inverting}.
According to a theorem of Koch and Venkov
\cite{KoVe1975},
the Galois group \(\mathrm{Gal}(\mathrm{F}_p^\infty(K)/K)\)
of the maximal unramified pro-\(p\) extension of an
imaginary quadratic field \(K=\mathbb{Q}(\sqrt{d})\)
must be a Schur \(\sigma\)-group
when \(p\) is an odd prime number.
Various deeper results on Schur \(\sigma\)-groups were established by Arrigoni
\cite{Ag1998}.

The focus of the present article is on the smallest odd prime \(p=3\)
and finite \(3\)-groups \(G\) with
elementary bicyclic commutator quotient \(G/G^\prime\simeq C_3\times C_3\)
and transfer kernel type \(\mathrm{F}\),
according to the terminology of Scholz and Taussky
\cite[\S\ 2.I.F, p. 36]{SoTa1934}.
Since the \(3\)-group \(G\) has two generators \(G=\langle x,y\rangle\)
with \(x^3,y^3\in G^\prime\),
it possesses four maximal self-conjugate subgroups,
\begin{equation}
\label{eqn:Maximal}
U_1=\langle x,G^\prime\rangle,\ U_2=\langle y,G^\prime\rangle,\ U_3=\langle xy,G^\prime\rangle,\ U_4=\langle xy^2,G^\prime\rangle,
\end{equation}
with Artin transfer homomorphisms \(T_i:\,G/G^\prime\to U_i/U_i^\prime\),
and the \textit{transfer kernel type} of \(G\) is defined by
\(\varkappa(G)=(\varkappa_1,\ldots,\varkappa_4)\), where
\(\varkappa_i=0\) if \(\ker(T_i)=G/G^\prime\) and
\(\varkappa_i=j\) if \(\ker(T_i)=U_j/G^\prime\) with \(1\le i,j\le 4\).
More precisely, \(\varkappa(G)\) is an orbit under the action of the symmetric group \(S_4\),
consisting of equivalent quartets \(\varkappa\sim\pi^{-1}\varkappa\pi\) with \(\pi\in S_4\)
\cite[\S\ 2.2, pp. 475--476]{Ma2012}.
Type \(\mathrm{F}\) splits into four subtypes,
\(\mathrm{F}.7\), \(\varkappa\sim (3443)\),
\(\mathrm{F}.11\), \(\varkappa\sim (1143)\),
\(\mathrm{F}.12\), \(\varkappa\sim (1343)\),
\(\mathrm{F}.13\), \(\varkappa\sim (3143)\),
sharing a common transposition \((43)\).
Types \(\mathrm{F}.11\), \(\mathrm{F}.12\) have a fixed point \((1)\),
types \(\mathrm{F}.7\), \(\mathrm{F}.13\) do not.

The smallest metabelian \(\sigma\)-groups \(M\) of type \(\mathrm{F}\) have
order \(\mathrm{ord}(M)=3^{9}\), coclass \(\mathrm{cc}(M)=4\) and nilpotency class \(\mathrm{cl}(M)=5\).
The smallest non-metabelian Schur \(\sigma\)-groups \(S\) of type \(\mathrm{F}\) have
order \(\mathrm{ord}(S)=3^{20}\), soluble length \(\mathrm{sl}(S)=3\),
coclass \(\mathrm{cc}(S)=11\) and nilpotency class \(\mathrm{cl}(S)=9\).
Since the order of \(3\)-groups \(G\) in the SmallGroups database
\cite{BEO2005}
is bounded by \(\mathrm{ord}(G)\le 3^{8}\),
all the groups under investigation must be constructed by the \(p\)-group generation algorithm
\cite{Nm1977,Ob1990,HEO2005},
and they will be denoted in the terminology of the ANUPQ-package
\cite{GNO2006}
by the combination of an absolute identifier and one (or several) relative identifier(s) in the shape
\begin{equation}
\label{eqn:Identifier}
% G\simeq\mathrm{SmallGroup}(o,i)-\#s;j \text{ with order } o\le 3^{8}, \text{ step size } s\le\nu, \text{ and } i,j\in\mathbb{N},
G\simeq\langle o,i\rangle-\#s;j \text{ with order } o\le 3^{8}, \text{ step size } s\le\nu, \text{ and } i,j\in\mathbb{N},
\end{equation}
where \(i\) is the absolute identifier in the SmallGroups library, and
\(\nu\) is the \textit{nuclear rank} of \(\langle o,i\rangle\).

We begin with basic foundations concerning
the evolution of transfer kernels in descendant trees in \S\
\ref{s:Evolution},
second abelian type invariants in \S\
\ref{s:ATI2},
and bifurcation in descendant trees in \S\
\ref{s:Bifurcation}.

In \S\
\ref{s:SporCc4},
the smallest \textit{sporadic} metabelian \(\sigma\)-groups \(M\) of type \(\mathrm{F}\)
\textit{outside of coclass trees} are identified.
They are all of order \(3^{9}\) and coclass \(4\).
By the construction of \textit{extremal root paths},
suitable non-metabelian Schur \(\sigma\)-groups \(S\) with \(S/S^{\prime\prime}\simeq M\) are provided.
The smallest of them are of order \(3^{20}\) and soluble length \(3\).

In \S\
\ref{s:PeriCc4},
the smallest \textit{periodic} metabelian \(\sigma\)-groups \(M\) of type \(\mathrm{F}\)
are identified as \textit{vertices of coclass trees}
\cite{Ma2015,Ma2016,Ma2018a}.
They are all of order \(3^{11}\) and coclass \(4\).
Smallest non-metabelian Schur \(\sigma\)-groups \(S\) with \(S/S^{\prime\prime}\simeq M\)
are of order \(3^{23}\) and soluble length \(3\).

In \S\
\ref{s:SporCc6},
bigger \textit{sporadic} metabelian \(\sigma\)-groups \(M\) of type \(\mathrm{F}\)
outside of coclass trees are identified.
They are all of order \(3^{13}\) and coclass \(6\).
Smallest non-metabelian Schur \(\sigma\)-groups \(S\) with \(S/S^{\prime\prime}\simeq M\)
are of order \(3^{26}\) and soluble length \(3\).

%\newpage
%--------------------------------------------------------------------------------

\section{Evolution (Contraction) of Transfer Kernels in Trees}
\label{s:Evolution}

\noindent
In
\cite{BuMa2015},
and with more depth in
\cite{Ma2015},
we have seen that finite \(3\)-groups
with transfer kernel types
\(\mathrm{E}.6\), \(\mathrm{E}.14\), respectively
\(\mathrm{E}.8\), \(\mathrm{E}.9\),
are immediate descendants of \textit{skeleton groups}
with transfer kernel type \(\mathrm{c}.18\),
respectively \(\mathrm{c}.21\). 
For finite \(3\)-groups
with one of the four transfer kernel types \(\mathrm{F}\),
we have to establish a corresponding overview of possible skeleton types.

\begin{lemma}
\label{lem:Skeleton}
\textbf{(Skeleton types.)}
Sporadic metabelian \(3\)-groups, as isolated vertices outside of coclass trees,
with one of the four transfer kernel types \(\mathrm{F}\)
are immediate descendants with step size \(s=2\)
of skeleton groups with transfer kernel type \(\mathrm{b}.10\),
\(\varkappa\sim (0043)\). \\
Periodic metabelian \(3\)-groups, as terminal vertices of coclass trees,
with one of the four transfer kernel types \(\mathrm{F}\)
are immediate descendants with step size \(s=1\)
of skeleton groups with three possible transfer kernel types \(\mathrm{d}\).
In detail:
\begin{itemize}
\item
Type \(\mathrm{F}.7\), \(\varkappa(K)\sim (4343)\),
can only arise from type \(\mathrm{d}.19\), \(\varkappa(K)\sim (0343)\).
\item
Type \(\mathrm{F}.11\), \(\varkappa(K)\sim (2243)\sim (1143)\), \\
arises either from type \(\mathrm{d}.23\), \(\varkappa(K)\sim (0243)\),
or from type \(\mathrm{d}.25\), \(\varkappa(K)\sim (0143)\),
\item
Type \(\mathrm{F}.12\), \(\varkappa(K)\sim (1343)\sim (3243)\), \\
arises either from type \(\mathrm{d}.19\), \(\varkappa(K)\sim (0343)\),
or from type \(\mathrm{d}.23\), \(\varkappa(K)\sim (0243)\),
\item
Type \(\mathrm{F}.13\), \(\varkappa(K)\sim (2343)\sim (3143)\), \\
arises either from type \(\mathrm{d}.19\), \(\varkappa(K)\sim (0343)\),
or from type \(\mathrm{d}.25\), \(\varkappa(K)\sim (0143)\).
\end{itemize}
\end{lemma}

\begin{proof}
By the \textit{antitony principle} for the Artin pattern
\cite[\S\S\ 5.1--5.4, pp. 78--87]{Ma2016},
a partial (one-dimensional) transfer kernel is inherited
by an immediate descendant from its parent,
whereas a total (two-dimensional) transfer kernel
may contract to any one-dimensional transfer kernel.

For a sporadic isolated vertex of step size \(s=2\),
both total kernels of type \(\mathrm{b}.10\),
\(\varkappa\sim (0043)\),
of the parent may shrink at once,
and thus the immediate descendant may have any
of the four transfer kernel types \(\mathrm{F}\),
\(\varkappa(K)\sim (4343)\) or \(\varkappa(K)\sim (1143)\) or \(\varkappa(K)\sim (1343)\) or \(\varkappa(K)\sim (2343)\).

For a periodic isolated vertex of step size \(s=1\),
the single total kernel of a parent of
type \(\mathrm{d}.19\), \(\varkappa(K)\sim (0343)\),
may either contract to
type \(\mathrm{F}.12\), \(\varkappa(K)\sim (1343)\), or
type \(\mathrm{F}.13\), \(\varkappa(K)\sim (2343)\) or
type \(\mathrm{H}.4\), \(\varkappa(K)\sim (3343)\).
The latter must be eliminated, when type \(\mathrm{F}\) is desired.

The single total kernel of a parent of
type \(\mathrm{d}.23\), \(\varkappa(K)\sim (0243)\),
may either shrink to
type \(\mathrm{F}.12\), \(\varkappa(K)\sim (3243)\), or
type \(\mathrm{F}.11\), \(\varkappa(K)\sim (2243)\) or
type \(\mathrm{G}.16\), \(\varkappa(K)\sim (1243)\).
The latter with two fixed points must be eliminated, when type \(\mathrm{F}\) is desired.

The single total kernel of a parent of
type \(\mathrm{d}.25\), \(\varkappa(K)\sim (0143)\),
may either shrink to
type \(\mathrm{F}.11\), \(\varkappa(K)\sim (1143)\), or
type \(\mathrm{F}.13\), \(\varkappa(K)\sim (3143)\) or
type \(\mathrm{G}.19\), \(\varkappa(K)\sim (2143)\).
The latter with two transpositions must be eliminated, when type \(\mathrm{F}\) is desired.
\end{proof}

%\newpage
%--------------------------------------------------------------------------------

\section{Second Abelian Type Invariants of Imaginary Quadratic Fields}
\label{s:ATI2}

%\noindent
Let \(N(\ell)\) be the \textit{nearly homocyclic abelian \(3\)-group}
of logarithmic order \(\ell\), that is,
\(N(2j)=(j,j)\) if \(\ell=2j\ge 0\) is even, and
\(N(2j+1)=(j+1,j)\) if \(\ell=2j+1\ge 1\) is odd.
\textit{Second abelian type invariants} (ATI) \(\alpha^{(2)}(k)\)
of any algebraic number field \(k\) with \(3\)-class group \(\mathrm{Cl}_3(k)\simeq C_3\times C_3\),
second \(3\)-class group \(M=\mathrm{G}_3^2(k)=\mathrm{Gal}(\mathrm{F}_3^2(k)/k)\),
nilpotency class \(c=\mathrm{cl}(M)\), coclass \(r=\mathrm{cc}(M)\),
and capitulation type \(\mathrm{F}\)
are always of the logarithmic shape
\begin{equation}
\alpha^{(2)}(k)=\lbrack 11;(Pol;Com,A_1),(Cpl;Com,A_2),(111;Com,A_3),(111;Com,A_4)\rbrack
\end{equation}
with triplets \(A_1\), \(A_2\) and dodecuplets \(A_3\), \(A_4\), and
abelian quotient invariants
\(Pol=N(c)\)
of the polarization,
\(Cpl=N(r+1)\)
of the co-polarization
\cite[Eqn. (5.1), p. 140]{Ma2017}, and
\(Com=N(c-1)\times N(r-1)\)
of the commutator subgroup
\cite[Appendix]{Ma2014}.
(The defect is zero for type \(\mathrm{F}\).)

%--------------------------------------------------------------------------------

\begin{definition}
\label{dfn:Categories}
Second abelian type invariants \(\alpha^{(2)}(k)\)
are said to be of the
\begin{itemize}
\item
\textit{first category} (extreme) if all components of the dodecuplets \(A_3\), \(A_4\)
have at least rank \(4\),
\item
\textit{second category} (wild, elevated) \\
if some but not all components of the dodecuplets \(A_3\), \(A_4\)
have at least rank \(4\),
\item
\textit{third category} (tame, moderate) if all components of the dodecuplets \(A_3\), \(A_4\)
have rank \(3\).
\end{itemize}
\end{definition}

%--------------------------------------------------------------------------------

\noindent
Recent investigations by Eric Ahlqvist and Magnus Carlsson
suggest that imaginary quadratic fields \(K=\mathbb{Q}(\sqrt{d})\), \(d<0\),
with \(\alpha^{(2)}(K)\) of the first and second category 
have infinite \(3\)-class field towers.
Unfortunately, their proof cannot be applied to elementary bicyclic \(\mathrm{Cl}_3(K)\).
So the infinitude of the \(3\)-class field tower
for wild and extreme type \(\mathrm{F}\) of \(K\) remains a conjecture.

%\newpage
%--------------------------------------------------------------------------------

\section{Bifurcation in Descendant Trees}
\label{s:Bifurcation}

\noindent
The basic genetics of transfer kernel types in \S\
\ref{s:Evolution}
can now be applied to describe the mutual position of a
finite Schur \(\sigma\)-group \(S\) of type \(\mathrm{F}\) and category \(3\)
and its metabelianization \(M=S/S^{\prime\prime}\).
As opposed to the types \(\mathrm{G}\) and \(\mathrm{H}\),
where \(S\) may be a descendant of \(M\)
\cite{Ma2017},
the situation for the types \(\mathrm{E}\) and \(\mathrm{F}\)
can be described in terms of \textit{root paths}
\cite[Dfn. 2.2, p. 87, \S\ 3, p. 88]{Ma2018c}
by the following definition.

\begin{definition}
\label{dfn:Fork}
The vertex \(B\) of the biggest order
in the meet of the root paths of \(S\) and \(M\)
is called the \textit{bifurcation} between \(S\) and \(M\).
The vertex \(R\) of the smallest order
in the root path of \(S\)
which has the transfer kernel type of \(S\)
is called the \textit{settled root} of \(S\)
(since no further contraction of transfer kernels
can occur for descendants of \(R\)).
The \textit{fork} between \(S\) and \(M\)
is the union of the path from \(R\) to \(B\)
and the path from \(M\) to \(B\)
(which intersect in \(B\)).
\end{definition}

\noindent
In ostensively simplified form,
the vertices of the fork are only labelled by
the letters of their transfer kernel types,
and the edges are labelled by their step sizes.
In view of possible arithmetical realizations,
we shall restrict our investigations to
three main situations,
\begin{itemize}
\item
Schur \(\sigma\)-group \(S\) with sporadic metabelianization \(M\) of coclass \(\mathrm{cc}(M)=4\) in \S\
\ref{s:SporCc4},
\item
Schur \(\sigma\)-group \(S\) with periodic metabelianization \(M\) of coclass \(\mathrm{cc}(M)=4\) in \S\
\ref{s:PeriCc4}, and
\item
Schur \(\sigma\)-group \(S\) with sporadic metabelianization \(M\) of coclass \(\mathrm{cc}(M)=6\) in \S\
\ref{s:SporCc6}.
\end{itemize}

%--------------------------------------------------------------------------------

\begin{theorem}
\label{thm:Fork}
Let \(S=\mathrm{Gal}(\mathrm{F}_3^\infty(K)/K)\) be the Schur \(\sigma\)-Galois group
of the maximal unramified pro-\(3\) extension of an
imaginary quadratic number field \(K=\mathbb{Q}(\sqrt{d})\)
with fundamental discriminant \(d<0\),
elementary bicyclic \(3\)-class group \(\mathrm{Cl}_3(K)\simeq C_3\times C_3\),
capitulation type \(\mathrm{F}\), and
second abelian type invariants \(\alpha^{(2)}(K)\) of category \(3\).
Then the fork between \(S\) and
\(M=S/S^{\prime\prime}\simeq\mathrm{Gal}(\mathrm{F}_3^2(K)/K)\),
in dependence on the first abelian type invariants \(\alpha^{(1)}(K)\),
is given by
\begin{equation}
\label{eqn:Fork}
\begin{aligned}
\mathrm{F}\buildrel 2\over\to\mathrm{b}\buildrel 4\over\leftarrow\mathrm{F} \quad
& \text{ if } \alpha^{(1)}(K)=\lbrack 11;32,32,111,111\rbrack, \\
\mathrm{F}\buildrel 1\over\to\mathrm{d}\buildrel 1\over\to\mathrm{d}\buildrel 2\over\to\mathrm{b}
\buildrel 4\over\leftarrow\mathrm{d}\buildrel 2\over\leftarrow\mathrm{d}\buildrel 4\over\leftarrow\mathrm{F} \quad
& \text{ if } \alpha^{(1)}(K)=\lbrack 11;43,32,111,111\rbrack, \\
\mathrm{F}\buildrel 2\over\to\mathrm{b}\buildrel 2\over\to\mathrm{b}\buildrel 2\over\to\mathrm{b}
\buildrel 4\over\leftarrow\mathrm{b}\buildrel 2\over\leftarrow\mathrm{b}\buildrel 4\over\leftarrow\mathrm{F} \quad
& \text{ if } \alpha^{(1)}(K)=\lbrack 11;43,43,111,111\rbrack.
\end{aligned}
\end{equation}
\end{theorem}

\begin{proof}
In
\cite[\S\ 3.6, pp. 101--104]{Ma2018c}
the fork was called \textit{fork topology}.
The metabelian part of the fork
has been described in
\cite[Cor. 3.1, p. 104]{Ma2018c}
for the two sporadic situations
\(\alpha^{(1)}(K)=\lbrack 11;32,32,111,111\rbrack\)
with class \(5\), coclass \(4\), and
\(\alpha^{(1)}(K)=\lbrack 11;43,43,111,111\rbrack\)
with class \(7\), coclass \(6\), and in
\cite[Thm. 3.6, p. 103]{Ma2018c}
for the periodic situation
\(\alpha^{(1)}(K)=\lbrack 11;43,32,111,111\rbrack\)
with class \(7\), coclass \(4\).
The bifurcation is always the fixed vertex \(B=\langle 2187,64\rangle\).
Therefore, the root paths \(P\) in
\cite[pp. 103--104]{Ma2018c}
must be shortened by the ommission of the trailing three edges
\(\langle 2187,64\rangle\buildrel 2\over\to\langle 243,3\rangle\buildrel 2\over\to\langle 27,3\rangle\buildrel 1\over\to\langle 9,2\rangle\)
of type
\(\mathrm{b}\buildrel 2\over\to\mathrm{b}\buildrel 2\over\to\mathrm{a}\buildrel 1\over\to\mathrm{a}\).
The non-metabelian part of the fork
is a consequence of the \textit{extremal root path property}
of finite Schur \(\sigma\)-groups \(S\) of category \(3\),
and it is symmetric with respect to the types
(but not with respect to the step sizes).
\end{proof}

%\newpage
%--------------------------------------------------------------------------------

\section{Realization of Sporadic Metabelian \(\sigma\)-Groups \(\mathrm{G}_3^2(K)\) of Coclass \(4\)}
\label{s:SporCc4}

\noindent
In Table
\ref{tbl:ATI2cc4Spor},
second abelian type invariants \(\alpha^{(2)}(K)\)
of imaginary quadratic fields \(K=\mathbb{Q}(\sqrt{d})\)
with discriminants in the narrow range \(-5\cdot 10^5<d<0\),
bicyclic \(3\)-class group \(\mathrm{Cl}_3(K)\simeq C_3\times C_3\),
second \(3\)-class group \(M=\mathrm{G}_3^2(K)\), \(\mathrm{ord}(M)=3^9\),
and type \(\mathrm{F}\) are given,
ordered by capitulation subtypes \(\varkappa(K)\).
\textbf{Boldface} font emphasizes exceptional components
of the first or second category.

%\newpage
%--------------------------------------------------------------------------------

\renewcommand{\arraystretch}{1.0}

\begin{table}[ht]
\caption{\(2^{\text{nd}}\) ATI of \(K=\mathbb{Q}(\sqrt{d})\), \(d<0\), with sporadic \(M=\mathrm{G}_3^2(K)\), \(\mathrm{cc}(M)=4\)}
\label{tbl:ATI2cc4Spor}
\begin{center}
\begin{tabular}{|r|c|c|c|c|}
\hline
         Type & \multicolumn{4}{|c|}{\(\alpha^{(2)}(K)=\lbrack 11;(32;2221,A_1),(32;2221,A_2),(111;2221,A_3),(111;2221,A_4)\rbrack\)} \\
       \(-d\) & \(A_1\) & \(A_2\) & \(A_3\) & \(A_4\) \\
\hline
 \(\mathrm{F}.7\) & \multicolumn{4}{|c|}{\(\varkappa(K)=(3443)\)} \\
\hline
 \(124\,363\) & \(\mathbf{(4211)^3}\)  & \(\mathbf{(3211)^3}\)  & \(\mathbf{(1^6)^3,(2221)^3,(2211)^6}\)  & \(\mathbf{(2^21^3)^3,(21^4)^3,(2211)^6}\) \\
 \(225\,299\) & \((3111)^3\)           & \((3111)^3\)           & \((221)^3,(211)^9\)                     & \((221)^3,(211)^9\)                       \\
 \(260\,515\) & \(\mathbf{(3321)^3}\)  & \(\mathbf{(3321)^3}\)  & \(\mathbf{(1^6)^3,(3221)^3,(2221)^3,(2211)^3}\) & \(\mathbf{(1^6)^3,(21^4)^3,(3221)^3,(2211)^3}\) \\
 \(343\,380\) & \((3111)^3\)           & \((3111)^3\)           & \((221)^3,(211)^9\)                     & \((221)^3,(211)^9\)                       \\
 \(423\,476\) & \((3111)^3\)           & \((3111)^3\)           & \((221)^3,(211)^9\)                     & \((221)^3,(211)^9\)                       \\
 \(486\,264\) & \((3111)^3\)           & \((3111)^3\)           & \((221)^3,(211)^9\)                     & \((221)^3,(211)^9\)                       \\
\hline
 \(\mathrm{F}.11\) & \multicolumn{4}{|c|}{\(\varkappa(K)=(1143)\)} \\
\hline
  \(27\,156\) & \((4111)^3\)           & \((3111)^3\)           & \((221)^3,(211)^9\)                     & \((221)^3,(211)^9\)                       \\
 \(241\,160\) & \((4111)^3\)           & \((3111)^3\)           & \((221)^3,(211)^9\)                     & \((221)^3,(211)^9\)                       \\
 \(394\,999\) & \((\mathbf{3}111)^3\)  & \((3111)^3\)           & \(\mathbf{(2111)^9},(221)^3\)           & \(\mathbf{(1^5)^3,(2111)^3,(1111)^6}\)    \\
 \(477\,192\) & \((4111)^3\)           & \((3111)^3\)           & \((221)^3,(211)^9\)                     & \((221)^3,(211)^9\)                       \\
 \(484\,804\) & \((4111)^3\)           & \((3111)^3\)           & \((221)^3,(211)^9\)                     & \((221)^3,(211)^9\)                       \\
\hline
 \(\mathrm{F}.12\) & \multicolumn{4}{|c|}{\(\varkappa(K)=(1343)\)} \\
\hline
  \(31\,908\) & \((3111)^3\)           & \((3111)^3\)           & \(\mathbf{(2111)^9},(221)^3\)           & \(\mathbf{(2111)^6,(1111)^6}\)            \\
 \(135\,587\) & \((3111)^3\)           & \((3111)^3\)           & \(\mathbf{(2211)^3,(2111)^6},(221)^3\)  & \(\mathbf{(1^5)^3,(2111)^3,(1111)^6}\)    \\
 \(160\,403\) & \(\mathbf{(3211)^3}\)  & \(\mathbf{(3211)^3}\)  & \(\mathbf{(21^4)^6,(2222)^3,(2211)^3}\) & \(\mathbf{(1^5)^6,(2211)^3,(2111)^3}\)    \\
 \(184\,132\) & \((\mathbf{4}111)^3\)  & \((3111)^3\)           & \(\mathbf{(2211)^3,(2111)^6},(221)^3\)  & \(\mathbf{(2211)^3,(2111)^3,(1111)^6}\)   \\
 \(189\,959\) & \((3111)^3\)           & \((3111)^3\)           & \(\mathbf{(2211)^3,(2111)^6},(221)^3\)  & \(\mathbf{(2111)^6,(1111)^6}\)            \\
 \(291\,220\) & \((3111)^3\)           & \((3111)^3\)           & \((221)^3,(211)^9\)                     & \((221)^3,(211)^9\)                       \\
 \(454\,631\) & \((3111)^3\)           & \((3111)^3\)           & \(\mathbf{(2111)^9},(221)^3\)           & \(\mathbf{(2111)^6,(1111)^6}\)            \\
 \(499\,159\) & \((3111)^3\)           & \((3111)^3\)           & \(\mathbf{(2111)^9},(221)^3\)           & \(\mathbf{(1^5)^3,(2111)^3,(1111)^6}\)    \\
\hline
 \(\mathrm{F}.13\) & \multicolumn{4}{|c|}{\(\varkappa(K)=(3143)\)} \\
\hline
  \(67\,480\) & \((4111)^3\)           & \((3111)^3\)           & \(\mathbf{(3211)^3,(2111)^6},(221)^3\)  & \(\mathbf{(2111)^6,(1111)^6}\)            \\
 \(104\,627\) & \((4111)^3\)           & \((3111)^3\)           & \(\mathbf{(2111)^9},(221)^3\)           & \(\mathbf{(1^5)^3,(2111)^3,(1111)^6}\)    \\
 \(167\,064\) & \((4111)^3\)           & \((3111)^3\)           & \((221)^3,(211)^9\)                     & \((221)^3,(211)^9\)                       \\
 \(224\,580\) & \(\mathbf{(3211)^3}\)  & \(\mathbf{(3211)^3}\)  & \(\mathbf{(21^4)^3,(1^5)^6,(2211)^3}\)  & \(\mathbf{(1^5)^3,(2211)^3,(2111)^6}\) \\
 \(287\,155\) & \((4111)^3\)           & \((3111)^3\)           & \(\mathbf{(2111)^9},(221)^3\)           & \(\mathbf{(2111)^6,(1111)^6}\)            \\
 \(296\,407\) & \((4111)^3\)           & \((3111)^3\)           & \((221)^3,(211)^9\)                     & \((221)^3,(211)^9\)                       \\
 \(317\,747\) & \((4111)^3\)           & \((3111)^3\)           & \((221)^3,(211)^9\)                     & \((221)^3,(211)^9\)                       \\
 \(344\,667\) & \((4111)^3\)           & \((3111)^3\)           & \(\mathbf{(2211)^3,(2111)^6},(221)^3\)  & \(\mathbf{(2211)^3,(2111)^3,(1111)^6}\)   \\
 \(401\,603\) & \((4111)^3\)           & \((3111)^3\)           & \((221)^3,(211)^9\)                     & \((221)^3,(211)^9\)                       \\
 \(426\,891\) & \((4111)^3\)           & \((3111)^3\)           & \(\mathbf{(2111)^9},(221)^3\)           & \(\mathbf{(2111)^6,(1111)^6}\)            \\
 \(487\,727\) & \((\mathbf{3}111)^3\)  & \((3111)^3\)           & \(\mathbf{(2111)^9},(221)^3\)           & \(\mathbf{(1^5)^3,(2111)^3,(1111)^6}\)    \\
\hline
\end{tabular}
\end{center}
\end{table}

%\newpage
%--------------------------------------------------------------------------------

\noindent
Throughout the sequel, we restrict our investigations to the tame situation of category \(3\).

%--------------------------------------------------------------------------------

\bigskip
\noindent
The second AQI,
\(\alpha^{(2)}(S)=\lbrack 11;(32;2221,A_1),(32;2221,A_2),(111;2221,A_3),(111;2221,A_4)\rbrack\),
in the Tables
\ref{tbl:CorrespondenceSporCC4F11},
\ref{tbl:CorrespondenceSporCC4F12},
\ref{tbl:CorrespondenceSporCC4F13},
\ref{tbl:CorrespondenceSporCC4F7}
contain triplets \(A_1\), \(A_2\) and dodecuplets \(A_3\), \(A_4\).
In terms of the \textit{fixed bifurcation} \(B=\langle 2187,64\rangle\),
the tables give a complete overview of all possible Schur \(\sigma\)-groups \(S\)
among the descendants of \(R=B-\#4;k\) as non-metabelian roots,
for any assigned metabelian \(\sigma\)-group
\(M=B-\#2;j\) with \(M=R/R^{\prime\prime}=S/S^{\prime\prime}\),
order \(\mathrm{ord}(M)=3^{9}\), nilpotency class \(\mathrm{cl}(M)=5\), coclass \(\mathrm{cc}(M)=4\),
and transfer kernel types \(\mathrm{F}.11\), \(\mathrm{F}.12\), \(\mathrm{F}.13\), \(\mathrm{F}.7\).
All second AQI are of category \(3\) (tame, moderate).

%--------------------------------------------------------------------------------

\renewcommand{\arraystretch}{1.0}

\begin{table}[ht]
\caption{Correspondence for sporadic \(M=S/S^{\prime\prime}\), \(\mathrm{cc}(M)=4\), of type \(\mathrm{F}.11\)}
\label{tbl:CorrespondenceSporCC4F11}
\begin{center}
\begin{tabular}{|c|c||c|c|c|c||c|c|c|c|}
\hline
   Type &  &  &  &  &  & \multicolumn{4}{|c|}{\(\alpha^{(2)}(S)\)} \\
 lo\((M)\)& \(j\) & \(k\) & lo\((S)\) & sl\((S)\) & \(\#\) & \(A_1\) & \(A_2\) & \(A_3\) & \(A_4\) \\
\hline
 \(\mathrm{F}.11\) &  &  &  &  &  & \multicolumn{4}{|c|}{\(\varkappa(S)=(1143)\)} \\
\hline
  \(9\) & \(36\) & \(127\) & \(20\) & \(3\) & \(81/81\) & \((4111)^3\) & \((3111)^3\) & \((221)^3,(211)^9\) & \((221)^3,(211)^9\) \\
   \(\) &   \(\) & \(144\) & \(20\) & \(3\) & \(27/27\) & \((4111)^3\) & \((3111)^3\) & \((221)^3,(211)^9\) & \((221)^3,(211)^9\) \\
   \(\) &   \(\) & \(172\) & \(20\) & \(3\) & \(81/81\) & \((4111)^3\) & \((3111)^3\) & \((221)^3,(211)^9\) & \((221)^3,(211)^9\) \\
   \(\) &   \(\) & \(180\) & \(20\) & \(3\) & \(27/27\) & \((4111)^3\) & \((3111)^3\) & \((221)^3,(211)^9\) & \((221)^3,(211)^9\) \\
  \(9\) & \(38\) & \(119\) & \(20\) & \(3\) & \(81/81\) & \((4111)^3\) & \((3111)^3\) & \((221)^3,(211)^9\) & \((221)^3,(211)^9\) \\
   \(\) &   \(\) & \(139\) & \(20\) & \(3\) & \(27/27\) & \((4111)^3\) & \((3111)^3\) & \((221)^3,(211)^9\) & \((221)^3,(211)^9\) \\
   \(\) &   \(\) & \(164\) & \(20\) & \(3\) & \(81/81\) & \((4111)^3\) & \((3111)^3\) & \((221)^3,(211)^9\) & \((221)^3,(211)^9\) \\
   \(\) &   \(\) & \(182\) & \(20\) & \(3\) & \(27/27\) & \((4111)^3\) & \((3111)^3\) & \((221)^3,(211)^9\) & \((221)^3,(211)^9\) \\
\hline
\end{tabular}
\end{center}
\end{table}

%--------------------------------------------------------------------------------

\begin{theorem}
\label{thm:F11Category3}
Let \(K=\mathbb{Q}(\sqrt{d})\) be an imaginary quadratic number field
with discriminant \(d<0\),
bicyclic \(3\)-class group \(\mathrm{Cl}_3(K)\simeq C_3\times C_3\),
second \(3\)-class group \(M=\mathrm{G}_3^2(K)\), \(\mathrm{ord}(M)=3^9\),
capitulation type \(\mathrm{F}.11\), \(\varkappa(K)\sim (1143)\),
and first abelian type invariants
\(\alpha^{(1)}(K)=\lbrack 11;32,32,111,111\rbrack\).
If the second abelian type invariants are of third category,
they must have the shape
\begin{equation}
\label{eqn:F11Category3}
\alpha^{(2)}(K)=\lbrack 11;(32;2221,(4111)^3),(32;2221,(3111)^3),(111;2221,(221)^3,(211)^9)^2\rbrack,
\end{equation}
and the \(3\)-class field tower \(\mathrm{F}_3^\infty(K)\) has precise length \(\ell_3(K)=3\) and relative degree \(3^{20}\) over \(K\).
\end{theorem}

\begin{proof}
This is an immediate consequence of Table
\ref{tbl:CorrespondenceSporCC4F11},
since all the non-metabelian roots \(R=B-\#4;k\)
with \(k\in\lbrace 119,127,139,144,164,172,180,182\rbrace\)
lead to Schur \(\sigma\)-groups \(S\) with soluble length \(\mathrm{sl}(S)=3\),
logarithmic order \(\mathrm{lo}(S)=20\),
and have the second abelian type invariants in Formula
\eqref{eqn:F11Category3}.
\end{proof}

%--------------------------------------------------------------------------------

\begin{example}
\label{exm:F11Category3}
According to Table
\ref{tbl:ATI2cc4Spor},
the quadratic fields \(K=\mathbb{Q}(\sqrt{d})\) with discriminants
\(d\in\lbrace -27\,156,-241\,160,-477\,192,-484\,804\rbrace\)
have the capitulation type \(\mathrm{F}.11\), \(\varkappa(K)=(1143)\),
and second abelian type invariants of third category in Formula
\eqref{eqn:F11Category3}.
According to Theorem
\ref{thm:F11Category3},
they have a \(3\)-class field tower \(\mathrm{F}_3^\infty(K)\) of precise length \(\ell_3(K)=3\).
\end{example}

%\newpage
%--------------------------------------------------------------------------------

\renewcommand{\arraystretch}{1.0}

\begin{table}[ht]
\caption{Correspondence for sporadic \(M=S/S^{\prime\prime}\), \(\mathrm{cc}(M)=4\), of type \(\mathrm{F}.12\)}
\label{tbl:CorrespondenceSporCC4F12}
\begin{center}
\begin{tabular}{|c|c||c|c|c|c||c|c|c|c|}
\hline
   Type &  &  &  &  &  & \multicolumn{4}{|c|}{\(\alpha^{(2)}(S)\)} \\
 lo\((M)\)& \(j\) & \(k\) & lo\((S)\) & sl\((S)\) & \(\#\) & \(A_1\) & \(A_2\) & \(A_3\) & \(A_4\) \\
\hline
 \(\mathrm{F}.12\) &  &  &  &  &  & \multicolumn{4}{|c|}{\(\varkappa(S)=(1343)\)} \\
\hline
  \(9\) & \(43\) & \(125\) & \(20\) & \(3\) & \(27/27\) & \((3111)^3\) & \((3111)^3\) & \((221)^3,(211)^9\) & \((221)^3,(211)^9\) \\
   \(\) &   \(\) & \(143\) & \(20\) & \(3\) & \(27/27\) & \((3111)^3\) & \((3111)^3\) & \((221)^3,(211)^9\) & \((221)^3,(211)^9\) \\
   \(\) &   \(\) & \(170\) & \(20\) & \(3\) & \(27/27\) & \((4111)^3\) & \((3111)^3\) & \((221)^3,(211)^9\) & \((221)^3,(211)^9\) \\
   \(\) &   \(\) & \(187\) & \(20\) & \(3\) & \(27/27\) & \((4111)^3\) & \((3111)^3\) & \((221)^3,(211)^9\) & \((221)^3,(211)^9\) \\
  \(9\) & \(46\) & \(130\) & \(20\) & \(3\) & \(27/27\) & \((3111)^3\) & \((3111)^3\) & \((221)^3,(211)^9\) & \((221)^3,(211)^9\) \\
   \(\) &   \(\) & \(146\) & \(20\) & \(3\) & \(27/27\) & \((3111)^3\) & \((3111)^3\) & \((221)^3,(211)^9\) & \((221)^3,(211)^9\) \\
   \(\) &   \(\) & \(175\) & \(20\) & \(3\) & \(27/27\) & \((4111)^3\) & \((3111)^3\) & \((221)^3,(211)^9\) & \((221)^3,(211)^9\) \\
   \(\) &   \(\) & \(190\) & \(20\) & \(3\) & \(27/27\) & \((4111)^3\) & \((3111)^3\) & \((221)^3,(211)^9\) & \((221)^3,(211)^9\) \\
  \(9\) & \(51\) & \(113\) & \(20\) & \(3\) & \(27/27\) & \((3111)^3\) & \((3111)^3\) & \((221)^3,(211)^9\) & \((221)^3,(211)^9\) \\
   \(\) &   \(\) & \(126\) & \(20\) & \(3\) & \(27/27\) & \((3111)^3\) & \((3111)^3\) & \((221)^3,(211)^9\) & \((221)^3,(211)^9\) \\
   \(\) &   \(\) & \(157\) & \(23,26,29\) & \(4\) & \(27/27\) & \((4111)^3\) & \((3111)^3\) & \((221)^3,(211)^9\) & \((221)^3,(211)^9\) \\
   \(\) &   \(\) & \(194\) & \(20\) & \(3\) & \(27/27\) & \((4111)^3\) & \((3111)^3\) & \((221)^3,(211)^9\) & \((221)^3,(211)^9\) \\
  \(9\) & \(53\) & \(116\) & \(20\) & \(3\) & \(27/27\) & \((3111)^3\) & \((3111)^3\) & \((221)^3,(211)^9\) & \((221)^3,(211)^9\) \\
   \(\) &   \(\) & \(118\) & \(20\) & \(3\) & \(27/27\) & \((3111)^3\) & \((3111)^3\) & \((221)^3,(211)^9\) & \((221)^3,(211)^9\) \\
   \(\) &   \(\) & \(160\) & \(23,26\) & \(4\) & \(27/27\) & \((4111)^3\) & \((3111)^3\) & \((221)^3,(211)^9\) & \((221)^3,(211)^9\) \\
   \(\) &   \(\) & \(195\) & \(20\) & \(3\) & \(27/27\) & \((4111)^3\) & \((3111)^3\) & \((221)^3,(211)^9\) & \((221)^3,(211)^9\) \\
\hline
\end{tabular}
\end{center}
\end{table}

%--------------------------------------------------------------------------------

\begin{theorem}
\label{thm:F12Category3}
An imaginary quadratic field \(K=\mathbb{Q}(\sqrt{d})\)
with fundamental discriminant \(d<0\),
elementary bicyclic \(3\)-class group \(\mathrm{Cl}_3(K)\),
second \(3\)-class group \(M=\mathrm{G}_3^2(K)\), \(\mathrm{ord}(M)=3^9\),
capitulation type \(\mathrm{F}.12\), \(\varkappa(K)\sim (1343)\),
and second abelian type invariants of third category
\begin{equation}
\label{eqn:F12Category3a}
\alpha^{(2)}(K)=\lbrack 11;(32;2221,(3111)^3),(32;2221,(3111)^3),(111;2221,(221)^3,(211)^9)^2\rbrack
\end{equation}
has a \(3\)-class field tower \(\mathrm{F}_3^\infty(K)\) of precise length \(\ell_3(K)=3\).
For second abelian type invariants
\begin{equation}
\label{eqn:F12Category3b}
\alpha^{(2)}(K)=\lbrack 11;(32;2221,(4111)^3),(32;2221,(3111)^3),(111;2221,(221)^3,(211)^9)^2\rbrack,
\end{equation}
the \(3\)-class field tower \(\mathrm{F}_3^\infty(K)\)
has two possible lengths \(\ell_3(K)\in\lbrace 3,4\rbrace\).
\end{theorem}

\begin{proof}
This is an immediate consequence of Table
\ref{tbl:CorrespondenceSporCC4F12},
since the unique non-metabelian roots \(R=B-\#4;k\) with \(k\in\lbrace 157,160\rbrace\)
which lead to Schur \(\sigma\)-groups \(S\) with soluble length \(\mathrm{sl}(S)=4\)
have the second abelian type invariants in Formula
\eqref{eqn:F12Category3b}.
\end{proof}

%--------------------------------------------------------------------------------

\begin{example}
\label{exm:F12Category3}
According to Table
\ref{tbl:ATI2cc4Spor},
the single quadratic field \(K=\mathbb{Q}(\sqrt{d})\) with discriminant
\(d=-291\,220\)
has capitulation type \(\mathrm{F}.12\), \(\varkappa(K)=(1343)\),
and second abelian type invariants of third category
\eqref{eqn:F12Category3a}.
According to Theorem
\ref{thm:F12Category3},
it has a \(3\)-class field tower \(\mathrm{F}_3^\infty(K)\) of precise length \(\ell_3(K)=3\).
\end{example}

%\newpage
%--------------------------------------------------------------------------------

\renewcommand{\arraystretch}{1.0}

\begin{table}[hb]
\caption{Correspondence for sporadic \(M=S/S^{\prime\prime}\), \(\mathrm{cc}(M)=4\), of type \(\mathrm{F}.13\)}
\label{tbl:CorrespondenceSporCC4F13}
\begin{center}
\begin{tabular}{|c|c||c|c|c|c||c|c|c|c|}
\hline
   Type &  &  &  &  &  & \multicolumn{4}{|c|}{\(\alpha^{(2)}(S)\)} \\
 lo\((M)\)& \(j\) & \(k\) & lo\((S)\) & sl\((S)\) & \(\#\) & \(A_1\) & \(A_2\) & \(A_3\) & \(A_4\) \\
\hline
 \(\mathrm{F}.13\) &  &  &  &  &  & \multicolumn{4}{|c|}{\(\varkappa(S)=(3143)\)} \\
\hline
  \(9\) & \(41\) & \(132\) & \(20\) & \(3\) & \(27/27\) & \((4111)^3\) & \((3111)^3\) & \((221)^3,(211)^9\) & \((221)^3,(211)^9\) \\
   \(\) &   \(\) & \(147\) & \(20\) & \(3\) & \(27/27\) & \((4111)^3\) & \((3111)^3\) & \((221)^3,(211)^9\) & \((221)^3,(211)^9\) \\
   \(\) &   \(\) & \(177\) & \(20\) & \(3\) & \(27/27\) & \((3111)^3\) & \((3111)^3\) & \((221)^3,(211)^9\) & \((221)^3,(211)^9\) \\
   \(\) &   \(\) & \(185\) & \(20\) & \(3\) & \(27/27\) & \((3111)^3\) & \((3111)^3\) & \((221)^3,(211)^9\) & \((221)^3,(211)^9\) \\
  \(9\) & \(47\) & \(122\) & \(20\) & \(3\) & \(27/27\) & \((4111)^3\) & \((3111)^3\) & \((221)^3,(211)^9\) & \((221)^3,(211)^9\) \\
   \(\) &   \(\) & \(141\) & \(20\) & \(3\) & \(27/27\) & \((4111)^3\) & \((3111)^3\) & \((221)^3,(211)^9\) & \((221)^3,(211)^9\) \\
   \(\) &   \(\) & \(167\) & \(20\) & \(3\) & \(27/27\) & \((3111)^3\) & \((3111)^3\) & \((221)^3,(211)^9\) & \((221)^3,(211)^9\) \\
   \(\) &   \(\) & \(192\) & \(20\) & \(3\) & \(27/27\) & \((3111)^3\) & \((3111)^3\) & \((221)^3,(211)^9\) & \((221)^3,(211)^9\) \\
  \(9\) & \(50\) & \(112\) & \(23,26,29\) & \(4\) & \(27/27\) & \((4111)^3\) & \((3111)^3\) & \((221)^3,(211)^9\) & \((221)^3,(211)^9\) \\
   \(\) &   \(\) & \(135\) & \(20\) & \(3\) & \(27/27\) & \((4111)^3\) & \((3111)^3\) & \((221)^3,(211)^9\) & \((221)^3,(211)^9\) \\
   \(\) &   \(\) & \(158\) & \(20\) & \(3\) & \(27/27\) & \((3111)^3\) & \((3111)^3\) & \((221)^3,(211)^9\) & \((221)^3,(211)^9\) \\
   \(\) &   \(\) & \(171\) & \(20\) & \(3\) & \(27/27\) & \((3111)^3\) & \((3111)^3\) & \((221)^3,(211)^9\) & \((221)^3,(211)^9\) \\
  \(9\) & \(52\) & \(115\) & \(23,26\) & \(4\) & \(27/27\) & \((4111)^3\) & \((3111)^3\) & \((221)^3,(211)^9\) & \((221)^3,(211)^9\) \\
   \(\) &   \(\) & \(137\) & \(20\) & \(3\) & \(27/27\) & \((4111)^3\) & \((3111)^3\) & \((221)^3,(211)^9\) & \((221)^3,(211)^9\) \\
   \(\) &   \(\) & \(161\) & \(20\) & \(3\) & \(27/27\) & \((3111)^3\) & \((3111)^3\) & \((221)^3,(211)^9\) & \((221)^3,(211)^9\) \\
   \(\) &   \(\) & \(163\) & \(20\) & \(3\) & \(27/27\) & \((3111)^3\) & \((3111)^3\) & \((221)^3,(211)^9\) & \((221)^3,(211)^9\) \\
\hline
\end{tabular}
\end{center}
\end{table}

%--------------------------------------------------------------------------------

\begin{theorem}
\label{thm:F13Category3}
An imaginary quadratic field \(K=\mathbb{Q}(\sqrt{d})\)
with fundamental discriminant \(d<0\),
elementary bicyclic \(3\)-class group \(\mathrm{Cl}_3(K)\),
second \(3\)-class group \(M=\mathrm{G}_3^2(K)\), \(\mathrm{ord}(M)=3^9\),
capitulation type \(\mathrm{F}.13\), \(\varkappa(K)\sim (3143)\),
and second abelian type invariants of third category
\begin{equation}
\label{eqn:F13Category3a}
\alpha^{(2)}(K)=\lbrack 11;(32;2221,(3111)^3),(32;2221,(3111)^3),(111;2221,(221)^3,(211)^9)^2\rbrack
\end{equation}
has a \(3\)-class field tower \(\mathrm{F}_3^\infty(K)\) of precise length \(\ell_3(K)=3\).
For second abelian type invariants
\begin{equation}
\label{eqn:F13Category3b}
\alpha^{(2)}(K)=\lbrack 11;(32;2221,(4111)^3),(32;2221,(3111)^3),(111;2221,(221)^3,(211)^9)^2\rbrack,
\end{equation}
the \(3\)-class field tower \(\mathrm{F}_3^\infty(K)\)
has two possible lengths \(\ell_3(K)\in\lbrace 3,4\rbrace\).
\end{theorem}

\begin{proof}
This is an immediate consequence of Table
\ref{tbl:CorrespondenceSporCC4F13},
since the unique non-metabelian roots \(R=B-\#4;k\) with \(k\in\lbrace 112,115\rbrace\)
which lead to Schur \(\sigma\)-groups \(S\) with soluble length \(\mathrm{sl}(S)=4\)
have the second abelian type invariants in Formula
\eqref{eqn:F13Category3b}.
\end{proof}

%--------------------------------------------------------------------------------

\begin{example}
\label{exm:F13Category3}
According to Table
\ref{tbl:ATI2cc4Spor},
the quadratic fields \(K=\mathbb{Q}(\sqrt{d})\) with discriminants
\(d\in\lbrace -167\,064,-296\,407,-317\,747,-401\,603\rbrace\)
have the capitulation type \(\mathrm{F}.13\), \(\varkappa(K)=(3143)\),
and second abelian type invariants of third category
\eqref{eqn:F13Category3b}.
According to Theorem
\ref{thm:F13Category3},
no sharp decision about the length of the \(3\)-class field tower \(\mathrm{F}_3^\infty(K)\)
can be drawn, since there are two possibilities \(\ell_3(K)\in\lbrace 3,4\rbrace\).
However, outside of the range in Table
\ref{tbl:ATI2cc4Spor},
we found two quadratic fields \(K=\mathbb{Q}(\sqrt{d})\) with discriminants
\(d\in\lbrace -731\,867,-803\,591\rbrace\),
capitulation type \(\mathrm{F}.13\), \(\varkappa(K)=(3143)\),
and second abelian type invariants of third category
\eqref{eqn:PeriF13Category3a}.
According to Theorem
\ref{thm:F13Category3},
they have a \(3\)-class field tower \(\mathrm{F}_3^\infty(K)\) of precise length \(\ell_3(K)=3\).
\end{example}

%\newpage
%--------------------------------------------------------------------------------

\renewcommand{\arraystretch}{1.0}

\begin{table}[ht]
\caption{Correspondence for sporadic \(M=S/S^{\prime\prime}\), \(\mathrm{cc}(M)=4\), of type \(\mathrm{F}.7\)}
\label{tbl:CorrespondenceSporCC4F7}
\begin{center}
\begin{tabular}{|c|c||c|c|c|c||c|c|c|c|}
\hline
   Type &  &  &  &  &  & \multicolumn{4}{|c|}{\(\alpha^{(2)}(S)\)} \\
 lo\((M)\)& \(j\) & \(k\) & lo\((S)\) & sl\((S)\) & \(\#\) & \(A_1\) & \(A_2\) & \(A_3\) & \(A_4\) \\
\hline
 \(\mathrm{F}.7\) &  &  &  &  &  & \multicolumn{4}{|c|}{\(\varkappa(S)=(3443)\)} \\
\hline
  \(9\) & \(55\) & \(121\) & \(20\) & \(3\) & \(27/27\) & \((3111)^3\) & \((3111)^3\) & \((221)^3,(211)^9\) & \((221)^3,(211)^9\) \\
   \(\) &   \(\) & \(131\) & \(20\) & \(3\) & \(27/27\) & \((3111)^3\) & \((3111)^3\) & \((221)^3,(211)^9\) & \((221)^3,(211)^9\) \\
   \(\) &   \(\) & \(165\) & \(26,29,32,38\) & \(4\) & \(27/27\) & \((3111)^3\) & \((3111)^3\) & \((221)^3,(211)^9\) & \((221)^3,(211)^9\) \\
   \(\) &   \(\) & \(196\) & \(20\) & \(3\) & \(27/27\) & \((3111)^3\) & \((3111)^3\) & \((221)^3,(211)^9\) & \((221)^3,(211)^9\) \\
  \(9\) & \(56\) & \(123\) & \(26,29,32,35,41\) & \(4\) & \(18/18\) & \((3111)^3\) & \((3111)^3\) & \((221)^3,(211)^9\) & \((221)^3,(211)^9\) \\
   \(\) &   \(\) & \(142\) & \(20\) & \(3\) & \(18/18\) & \((3111)^3\) & \((3111)^3\) & \((221)^3,(211)^9\) & \((221)^3,(211)^9\) \\
   \(\) &   \(\) & \(169\) & \(20\) & \(3\) & \(27/27\) & \((3111)^3\) & \((3111)^3\) & \((221)^3,(211)^9\) & \((221)^3,(211)^9\) \\
  \(9\) & \(58\) & \(128\) & \(26,29,32\) & \(4\) & \(18/18\) & \((3111)^3\) & \((3111)^3\) & \((221)^3,(211)^9\) & \((221)^3,(211)^9\) \\
   \(\) &   \(\) & \(145\) & \(20\) & \(3\) & \(18/18\) & \((3111)^3\) & \((3111)^3\) & \((221)^3,(211)^9\) & \((221)^3,(211)^9\) \\
   \(\) &   \(\) & \(174\) & \(20\) & \(3\) & \(27/27\) & \((3111)^3\) & \((3111)^3\) & \((221)^3,(211)^9\) & \((221)^3,(211)^9\) \\
\hline
\end{tabular}
\end{center}
\end{table}

%--------------------------------------------------------------------------------

\begin{theorem}
\label{thm:F7Category3}
For an imaginary quadratic field \(K=\mathbb{Q}(\sqrt{d})\)
with discriminant \(d<0\),
elementary bicyclic \(3\)-class group \(\mathrm{Cl}_3(K)\),
second \(3\)-class group \(M=\mathrm{G}_3^2(K)\), \(\mathrm{ord}(M)=3^9\),
and capitulation type \(\mathrm{F}.7\), \(\varkappa(K)\sim (3443)\),
the second abelian type invariants of third category
must have the shape
\begin{equation}
\label{eqn:F7Category3}
\alpha^{(2)}(K)=\lbrack 11;(32;2221,(3111)^3),(32;2221,(3111)^3),(111;2221,(221)^3,(211)^9)^2\rbrack,
\end{equation}
and the \(3\)-class field tower \(\mathrm{F}_3^\infty(K)\)
has two possible lengths \(\ell_3(K)\in\lbrace 3,4\rbrace\).
\end{theorem}

\begin{proof}
This is an immediate consequence of Table
\ref{tbl:CorrespondenceSporCC4F7},
since the unique non-metabelian roots \(R=B-\#4;k\) with \(k\in\lbrace 123,128,165\rbrace\)
which lead to Schur \(\sigma\)-groups \(S\) with soluble length \(\mathrm{sl}(S)=4\)
also have the common second abelian type invariants in Formula
\eqref{eqn:F7Category3}.
\end{proof}

%--------------------------------------------------------------------------------

\begin{example}
\label{exm:F7Category3}
According to Table
\ref{tbl:ATI2cc4Spor},
the quadratic fields \(K=\mathbb{Q}(\sqrt{d})\) with discriminants
\(d\in\lbrace -225\,299,-343\,380,-423\,476,-486\,264\rbrace\)
have the capitulation type \(\mathrm{F}.7\), \(\varkappa(K)=(3443)\),
and second abelian type invariants of third category
\eqref{eqn:F7Category3}.
According to Theorem
\ref{thm:F7Category3},
no sharp decision about the length of the \(3\)-class field tower \(\mathrm{F}_3^\infty(K)\)
can be drawn, since there are two possibilities \(\ell_3(K)\in\lbrace 3,4\rbrace\).
\end{example}

%--------------------------------------------------------------------------------

\begin{corollary}
\label{cor:SporadicCc4Category3}
The logarithmic order
of the Schur \(\sigma\)-group \(S=\mathrm{Gal}(\mathrm{F}_3^\infty(K)/K)\)
is generally \(\mathrm{lo}(S)=20\) for a \(3\)-class field tower of length \(\ell_3(K)=3\),
and it depends on the capitulation type \(\varkappa(K)\),
\[
\mathrm{lo}(S)\in
\begin{cases}
\lbrace 23,26,29\rbrace           & \text{ for type } \mathrm{F}.12 \text{ or } \mathrm{F}.13, \\
\lbrace 26,29,32,35,38,41\rbrace  & \text{ for type } \mathrm{F}.7,
\end{cases}
\] 
in the case of a tower of length \(\ell_3(K)=4\).
\end{corollary}

\begin{proof}
This is a common consequence of the Tables
\ref{tbl:CorrespondenceSporCC4F11},
\ref{tbl:CorrespondenceSporCC4F12},
\ref{tbl:CorrespondenceSporCC4F13},
\ref{tbl:CorrespondenceSporCC4F7}.
\end{proof}

%\newpage
%--------------------------------------------------------------------------------

\section{Realization of Periodic Metabelian \(\sigma\)-Groups \(\mathrm{G}_3^2(K)\) of Coclass \(4\)}
\label{s:PeriCc4}

\noindent
Table
\ref{tbl:ATI2cc4Peri}
gives second abelian type invariants \(\alpha^{(2)}(K)\) of imaginary quadratic fields \(K=\mathbb{Q}(\sqrt{d})\)
with discriminants in the bigger range \(-10^6<d<0\),
bicyclic \(3\)-class group \(\mathrm{Cl}_3(K)\simeq C_3\times C_3\),
second \(3\)-class group \(M=\mathrm{G}_3^2(K)\), \(\mathrm{ord}(M)=3^{11}\),
and type \(\mathrm{F}\),
ordered by capitulation subtypes \(\varkappa(K)\).
\textbf{Boldface} font indicates categories \(1\) and \(2\).

%\newpage
%--------------------------------------------------------------------------------

\renewcommand{\arraystretch}{1.0}

\begin{table}[ht]
\caption{\(2^{\text{nd}}\) ATI of \(K=\mathbb{Q}(\sqrt{d})\), \(d<0\), with periodic \(M=\mathrm{G}_3^2(K)\), \(\mathrm{cc}(M)=4\)}
\label{tbl:ATI2cc4Peri}
\begin{center}
\begin{tabular}{|r|c|c|c|c|}
\hline
         Type & \multicolumn{4}{|c|}{\(\alpha^{(2)}(K)=\lbrack 11;(43;3321,A_1),(32;3321,A_2),(111;3321,A_3),(111;3321,A_4)\rbrack\)} \\
       \(-d\) & \(A_1\) & \(A_2\) & \(A_3\) & \(A_4\) \\
\hline
 \(\mathrm{F}.7\) & \multicolumn{4}{|c|}{\(\varkappa(K)=(3443)\)} \\
\hline
 \(469\,816\) & \((\mathbf{5}211)^3\)  & \((3111)^3\)           & \(\mathbf{(2211)^3,(2111)^6},(221)^3\)    & \(\mathbf{(2211)^3,(2111)^3,(1111)^6}\)   \\
 \(643\,011\) & \((4211)^3\)           & \((3111)^3\)           & \((221)^3,(211)^9\)                       & \((221)^3,(211)^9\)                       \\
 \(797\,556\) & \((4211)^3\)           & \((3111)^3\)           & \((221)^3,(211)^9\)                       & \((221)^3,(211)^9\)                       \\
\hline
 \(\mathrm{F}.11\) & \multicolumn{4}{|c|}{\(\varkappa(K)=(1143)\)} \\
\hline
 \(469\,787\) & \((5211)^3\)           & \((3111)^3\)           & \((221)^3,(211)^9\)                       & \((221)^3,(211)^9\)                       \\
\hline
 \(\mathrm{F}.12\) & \multicolumn{4}{|c|}{\(\varkappa(K)\sim (1343)\)} \\
\hline
 \(249\,371\) & \(\mathbf{(4221)^3}\)  & \(\mathbf{(3211)^3}\)  & \(\mathbf{(21^4)^3,(2211)^3,(2111)^6}\)   & \(\mathbf{(21^4)^9,(2211)^3}\)            \\
 \(278\,427\) & \((5211)^3\)           & \((3111)^3\)           & \((221)^3,(211)^9\)                       & \((221)^3,(211)^9\)                       \\
 \(382\,123\) & \((4211)^3\)           & \((3111)^3\)           & \((221)^3,(211)^9\)                       & \((221)^3,(211)^9\)                       \\
\hline
 \(\mathrm{F}.13\) & \multicolumn{4}{|c|}{\(\varkappa(K)\sim (3143)\)} \\
\hline
 \(159\,208\) & \((4211)^3\)           & \((3111)^3\)           & \((221)^3,(211)^9\)                       & \((221)^3,(211)^9\)                       \\
 \(262\,628\) & \((5211)^3\)           & \((3111)^3\)           & \((221)^3,(211)^9\)                       & \((221)^3,(211)^9\)                       \\
 \(273\,284\) & \((4211)^3\)           & \((4111)^3\)           & \((221)^3,(211)^9\)                       & \((221)^3,(211)^9\)                       \\
 \(551\,112\) & \((5211)^3\)           & \((3111)^3\)           & \((221)^3,(211)^9\)                       & \((221)^3,(211)^9\)                       \\
 \(940\,943\) & \((4211)^3\)           & \((3111)^3\)           & \(\mathbf{(2111)^9},(221)^3\)             & \(\mathbf{(2111)^6,(1111)^6}\)            \\
 \(947\,463\) & \((4211)^3\)           & \((\mathbf{4}111)^3\)  & \(\mathbf{(2111)^9},(221)^3\)             & \(\mathbf{(2111)^6,(1111)^6}\)            \\
\hline
\end{tabular}
\end{center}
\end{table}

%\newpage
%--------------------------------------------------------------------------------

\noindent
Throughout the sequel, we restrict our investigations to the tame situation of category \(3\).

%--------------------------------------------------------------------------------

\bigskip
\noindent
For periodic metabelianizations \(M\) on coclass trees,
we need more space for relative identifiers in chains of descendants.
Since all second AQI to be considered for finite soluble length
are of third category (tame, moderate),
the second abelian quotient invariants,
\[
\alpha^{(2)}(S)=\lbrack 11;(43;3321,A_1),(32;3321,A_2),(111;3321,A_3),(111;3321,A_4)\rbrack,
\]
in the Tables
\ref{tbl:CorrespondencePeriCC4F13},
\ref{tbl:CorrespondencePeriCC4F12},
\ref{tbl:CorrespondencePeriCC4F11},
\ref{tbl:CorrespondencePeriCC4F7}
contain triplets \(A_1\), \(A_2\) and fixed dodecuplets \(A_3=A_4=\lbrack (221)^3,(211)^9\rbrack\),
excluded from the table in order to provide more space.
The tables give a complete overview of all possible Schur \(\sigma\)-groups \(S\)
among the descendants of \(R=B-\#4;k_1-\#2;k_2-\#4;k_3\) as non-metabelian roots,
for any assigned metabelian \(\sigma\)-group
\(M=B-\#2;j_1-\#1;j_2-\#1;j_3\)
with \(M=R/R^{\prime\prime}=S/S^{\prime\prime}\),
order \(\mathrm{ord}(M)=3^{11}\), nilpotency class \(\mathrm{cl}(M)=7\), coclass \(\mathrm{cc}(M)=4\),
and transfer kernel type \(\mathrm{F}.13\), \(\mathrm{F}.12\), \(\mathrm{F}.11\), \(\mathrm{F}.7\).

%--------------------------------------------------------------------------------

\begin{theorem}
\label{thm:PeriF13Category3}
For an imaginary quadratic field \(K=\mathbb{Q}(\sqrt{d})\)
with discriminant \(d<0\),
bicyclic \(3\)-class group \(\mathrm{Cl}_3(K)\simeq C_3\times C_3\),
second \(3\)-class group \(M=\mathrm{G}_3^2(K)\), \(\mathrm{ord}(M)=3^{11}\),
capitulation type \(\mathrm{F}.13\), \(\varkappa(K)\sim (3143)\),
and any among the second abelian type invariants of third category,
\begin{equation}
\label{eqn:PeriF13Category3a}
\alpha^{(2)}(K)=\lbrack 11;(43;3321,(4211)^3),(32;3321,(3111)^3),(111;3321,(221)^3,(211)^9)^2\rbrack,
\end{equation}
\begin{equation}
\label{eqn:PeriF13Category3b}
\alpha^{(2)}(K)=\lbrack 11;(43;3321,(5211)^3),(32;3321,(3111)^3),(111;3321,(221)^3,(211)^9)^2\rbrack,
\end{equation}
\begin{equation}
\label{eqn:PeriF13Category3c}
\alpha^{(2)}(K)=\lbrack 11;(43;3321,(4211)^3),(32;3321,(4111)^3),(111;3321,(221)^3,(211)^9)^2\rbrack,
\end{equation}
the \(3\)-class field tower \(\mathrm{F}_3^\infty(K)\)
has two possible lengths \(\ell_3(K)\in\lbrace 3,4\rbrace\).
\end{theorem}

\begin{proof}
This is an immediate consequence of Table
\ref{tbl:CorrespondencePeriCC4F13},
since for each of three patterns in Formulas
\eqref{eqn:PeriF13Category3a},
\eqref{eqn:PeriF13Category3b},
\eqref{eqn:PeriF13Category3c}
there exist Schur \(\sigma\)-groups \(S\) with soluble lengths \(\mathrm{sl}(S)\in\lbrace 3,4\rbrace\).
\end{proof}

%--------------------------------------------------------------------------------

\begin{example}
\label{exm:PeriF13Category3}
According to Table
\ref{tbl:ATI2cc4Peri},
the quadratic fields \(K=\mathbb{Q}(\sqrt{d})\) with discriminants
\(d\in\lbrace -159\,208,-262\,628,-273\,284,-551\,112\rbrace\)
have the capitulation type \(\mathrm{F}.13\), \(\varkappa(K)=(3143)\),
and \(2^{\text{nd}}\) ATI of third category
\eqref{eqn:PeriF13Category3a} for \(d=-159\,208\),
\eqref{eqn:PeriF13Category3b} for \(d\in\lbrace -262\,628,-551\,112\rbrace\),
\eqref{eqn:PeriF13Category3c} for \(d=-273\,284\).
According to Theorem
\ref{thm:PeriF13Category3},
no sharp decision about the length of the \(3\)-class field tower \(\mathrm{F}_3^\infty(K)\)
can be drawn, since there are two possibilities \(\ell_3(K)\in\lbrace 3,4\rbrace\).
However, disclosure of the skeleton type
\(\mathrm{d}.19\), \(j_1\in\lbrace 39,44\rbrace\), is enabled for
\eqref{eqn:PeriF13Category3b},
and \(\mathrm{d}.25\), \(j_1\in\lbrace 57,59\rbrace\), for
\eqref{eqn:PeriF13Category3c}.
\end{example}

%\newpage
%--------------------------------------------------------------------------------

\renewcommand{\arraystretch}{1.0}

\begin{table}[ht]
\caption{Correspondence for periodic \(M=S/S^{\prime\prime}\), \(\mathrm{cc}(M)=4\), of type \(\mathrm{F}.13\)}
\label{tbl:CorrespondencePeriCC4F13}
\begin{center}
\begin{tabular}{|c|c|c|c|c||c|c|c|c|c||c|c|}
\hline
     Type &  &  &  &  &  &  &  &  &  & \multicolumn{2}{|c|}{\(\alpha^{(2)}(S)\)} \\
lo\((M)\) & \(j_1\) & \(j_2\) & \(j_3\) & \(\#\) & \(k_1\) & \(k_3\) & lo\((S)\) & sl\((S)\) & \(\#\) & \(A_1\) & \(A_2\) \\
\hline
 \(\mathrm{F}.13\) &  &  &  &  &  &  &  &  &  & \multicolumn{2}{|c|}{\(\varkappa(S)=(3143)\)} \\
\hline
 \(11\) & \(39\) &  \(7\) &   \(3\) & \(1/25\) & \(117\) & \(24;19,26\) & \(29;32,41\) & \(3;4\) &  \(3/27\) & \((5211)^3\) & \((3111)^3\) \\
 \(11\) &   \(\) &   \(\) &    \(\) &     \(\) & \(138\) & \(11,15,16\) & \(23\)       & \(3\)   &  \(3/27\) & \((5211)^3\) & \((3111)^3\) \\
 \(11\) &   \(\) &   \(\) &    \(\) &     \(\) & \(162\) & \(12,13,17\) & \(26\)       & \(3\)   &  \(3/27\) & \((4211)^3\) & \((3111)^3\) \\
 \(11\) &   \(\) &   \(\) &    \(\) &     \(\) & \(184\) & \(19,24,26\) & \(23\)       & \(3\)   &  \(3/27\) & \((4211)^3\) & \((3111)^3\) \\

 \(11\) & \(39\) &  \(7\) &   \(8\) & \(1/25\) & \(117\) &    \(3;5,7\) & \(29;32\)    & \(3;4\) &  \(3/27\) & \((5211)^3\) & \((3111)^3\) \\
 \(11\) &   \(\) &   \(\) &    \(\) &     \(\) & \(138\) &    \(1,5,9\) & \(23\)       & \(3\)   &  \(3/27\) & \((5211)^3\) & \((3111)^3\) \\
 \(11\) &   \(\) &   \(\) &    \(\) &     \(\) & \(162\) &    \(2,6,7\) & \(26\)       & \(3\)   &  \(3/27\) & \((4211)^3\) & \((3111)^3\) \\
 \(11\) &   \(\) &   \(\) &    \(\) &     \(\) & \(184\) &    \(3,5,7\) & \(23\)       & \(3\)   &  \(3/27\) & \((4211)^3\) & \((3111)^3\) \\

 \(11\) & \(44\) &  \(1\) &   \(2\) & \(1/25\) & \(114\) &    \(6;2,7\) & \(29;32\)    & \(3;4\) &  \(3/27\) & \((5211)^3\) & \((3111)^3\) \\
 \(11\) &   \(\) &   \(\) &    \(\) &     \(\) & \(136\) &    \(3,5,7\) & \(23\)       & \(3\)   &  \(3/27\) & \((5211)^3\) & \((3111)^3\) \\
 \(11\) &   \(\) &   \(\) &    \(\) &     \(\) & \(159\) &    \(3,5,7\) & \(26\)       & \(3\)   &  \(3/27\) & \((4211)^3\) & \((3111)^3\) \\
 \(11\) &   \(\) &   \(\) &    \(\) &     \(\) & \(189\) & \(21,22,26\) & \(23\)       & \(3\)   &  \(3/27\) & \((4211)^3\) & \((3111)^3\) \\

 \(11\) & \(44\) &  \(1\) &   \(9\) & \(1/25\) & \(114\) & \(12;13,17\) & \(29;32,35\) & \(3;4\) &  \(3/27\) & \((5211)^3\) & \((3111)^3\) \\
 \(11\) &   \(\) &   \(\) &    \(\) &     \(\) & \(136\) & \(19,24,26\) & \(23\)       & \(3\)   &  \(3/27\) & \((5211)^3\) & \((3111)^3\) \\
 \(11\) &   \(\) &   \(\) &    \(\) &     \(\) & \(159\) & \(19,24,26\) & \(26\)       & \(3\)   &  \(3/27\) & \((4211)^3\) & \((3111)^3\) \\
 \(11\) &   \(\) &   \(\) &    \(\) &     \(\) & \(189\) &    \(1,5,9\) & \(23\)       & \(3\)   &  \(3/27\) & \((4211)^3\) & \((3111)^3\) \\
\hline
 \(11\) & \(57\) &  \(1\) &   \(2\) & \(1/15\) & \(124\) & \(10,15,17\) & \(26\)       & \(3\)   &  \(3/27\) & \((4211)^3\) & \((4111)^3\) \\
 \(11\) &   \(\) &   \(\) &    \(\) &     \(\) & \(124\) & \(20,22,27\) & \(23\)       & \(3\)   &  \(3/27\) & \((4211)^3\) & \((4111)^3\) \\
 \(11\) &   \(\) &   \(\) &    \(\) &     \(\) & \(168\) &    \(3,9;5\) & \(26;32\)    & \(4\)   &  \(3/15\) & \((4211)^3\) & \((3111)^3\) \\
 \(11\) &   \(\) &   \(\) &    \(\) &     \(\) & \(197\) &   \(5,9,12\) & \(29\)       & \(4\)   &  \(3/45\) & \((4211)^3\) & \((3111)^3\) \\
 \(11\) &   \(\) &   \(\) &    \(\) &     \(\) & \(197\) & \(16,22,25\) & \(29\)       & \(4\)   &  \(3/45\) & \((4211)^3\) & \((3111)^3\) \\
 \(11\) &   \(\) &   \(\) &    \(\) &     \(\) & \(197\) & \(28,32,42\) & \(29\)       & \(4\)   &  \(3/45\) & \((4211)^3\) & \((3111)^3\) \\ 
 
 \(11\) & \(57\) &  \(1\) &   \(4\) & \(1/15\) & \(124\) &    \(2,4,9\) & \(23\)       & \(3\)   &  \(3/27\) & \((4211)^3\) & \((4111)^3\) \\  
 \(11\) &   \(\) &   \(\) &    \(\) &     \(\) & \(124\) & \(12,14,16\) & \(26\)       & \(3\)   &  \(3/27\) & \((4211)^3\) & \((4111)^3\) \\
 \(11\) &   \(\) &   \(\) &    \(\) &     \(\) & \(168\) &  \(4,12;11\) & \(26;29\)    & \(4\)   &  \(3/15\) & \((4211)^3\) & \((3111)^3\) \\
 \(11\) &   \(\) &   \(\) &    \(\) &     \(\) & \(197\) &   \(2,7,17\) & \(29\)       & \(4\)   &  \(3/45\) & \((4211)^3\) & \((3111)^3\) \\
 \(11\) &   \(\) &   \(\) &    \(\) &     \(\) & \(197\) & \(19,23,29\) & \(29\)       & \(4\)   &  \(3/45\) & \((4211)^3\) & \((3111)^3\) \\
 \(11\) &   \(\) &   \(\) &    \(\) &     \(\) & \(197\) & \(33,37,41\) & \(29\)       & \(4\)   &  \(3/45\) & \((4211)^3\) & \((3111)^3\) \\
 
 \(11\) & \(59\) &  \(6\) &   \(3\) & \(1/15\) & \(129\) &    \(1,5,9\) & \(26\)       & \(3\)   &  \(3/27\) & \((4211)^3\) & \((4111)^3\) \\
 \(11\) &   \(\) &   \(\) &    \(\) &     \(\) & \(129\) & \(20,24,25\) & \(23\)       & \(3\)   &  \(3/27\) & \((4211)^3\) & \((4111)^3\) \\
 \(11\) &   \(\) &   \(\) &    \(\) &     \(\) & \(173\) &   \(6,14;2\) & \(26;29\)    & \(4\)   &  \(3/27\) & \((4211)^3\) & \((3111)^3\) \\
 \(11\) &   \(\) &   \(\) &    \(\) &     \(\) & \(198\) &   \(6,8,13\) & \(29\)       & \(4\)   &  \(3/45\) & \((4211)^3\) & \((3111)^3\) \\
 \(11\) &   \(\) &   \(\) &    \(\) &     \(\) & \(198\) & \(17,22,30\) & \(29\)       & \(4\)   &  \(3/45\) & \((4211)^3\) & \((3111)^3\) \\ 
 \(11\) &   \(\) &   \(\) &    \(\) &     \(\) & \(198\) & \(32,34,40\) & \(29\)       & \(4\)   &  \(3/45\) & \((4211)^3\) & \((3111)^3\) \\

 \(11\) & \(59\) &  \(6\) &   \(4\) & \(1/15\) & \(129\) &    \(3,8,4\) & \(29,32,35\) & \(4\)   &  \(3/27\) & \((4211)^3\) & \((4111)^3\) \\
 \(11\) &   \(\) &   \(\) &    \(\) &     \(\) & \(129\) & \(11,15,16\) & \(23\)       & \(3\)   &  \(3/27\) & \((4211)^3\) & \((4111)^3\) \\
 \(11\) &   \(\) &   \(\) &    \(\) &     \(\) & \(173\) &   \(4,12;8\) & \(26;29\)    & \(4\)   &  \(3/27\) & \((4211)^3\) & \((3111)^3\) \\
 \(11\) &   \(\) &   \(\) &    \(\) &     \(\) & \(198\) &   \(3,7,10\) & \(29\)       & \(4\)   &  \(3/45\) & \((4211)^3\) & \((3111)^3\) \\
 \(11\) &   \(\) &   \(\) &    \(\) &     \(\) & \(198\) & \(16,24,25\) & \(29\)       & \(4\)   &  \(3/45\) & \((4211)^3\) & \((3111)^3\) \\
 \(11\) &   \(\) &   \(\) &    \(\) &     \(\) & \(198\) & \(28,36,42\) & \(29\)       & \(4\)   &  \(3/45\) & \((4211)^3\) & \((3111)^3\) \\
\hline
\end{tabular}
\end{center}
\end{table}
%--------------------------------------------------------------------------------

\noindent
In the Tables
\ref{tbl:CorrespondencePeriCC4F13},
\ref{tbl:CorrespondencePeriCC4F12},
\ref{tbl:CorrespondencePeriCC4F11},
\ref{tbl:CorrespondencePeriCC4F7},
the Schur \(\sigma\)-groups \(S\)
have the settled root \(R=B-\#4;k_1-\#2;k_2-\#4;k_3\)
with some unique relative identifier \(k_2=k_2(k_1)\),
metabelianization \(M\simeq B-\#2;j_1-\#1;j_2-\#1;j_3\),
and \(2^{\text{nd}}\) AQI 
\(\alpha^{(2)}(S)=\lbrack 11;(43;3321,A_1),(32;3321,A_2),(111;3321,A_3),(111;3321,A_4)\rbrack\)
with \(A_3=A_4=((221)^3,(211)^9)\).
The relative identifier \(k_2\) cannot be given explicitly,
since the unique \(\sigma\)-descendant of \(B-\#4;k_1\)
is constructed directly by means of the \(p\)-covering group,
without using the \(p\)-group generation algorithm.

%\newpage
%--------------------------------------------------------------------------------

\renewcommand{\arraystretch}{1.0}

\begin{table}[ht]
\caption{Correspondence for periodic \(M=S/S^{\prime\prime}\), \(\mathrm{cc}(M)=4\), of type \(\mathrm{F}.12\)}
\label{tbl:CorrespondencePeriCC4F12}
\begin{center}
\begin{tabular}{|c|c|c|c|c||c|c|c|c|c||c|c|}
\hline
   Type &  &  &  &  &  &  &  &  &  & \multicolumn{2}{|c|}{\(\alpha^{(2)}(S)\)} \\
 lo\((M)\)& \(j_1\) & \(j_2\) & \(j_3\) & \(\#\) & \(k_1\) & \(k_3\) & lo\((S)\) & sl\((S)\) & \(\#\) & \(A_1\) & \(A_2\) \\
\hline
 \(\mathrm{F}.12\) &  &  &  &  &  &  &  &  &  & \multicolumn{2}{|c|}{\(\varkappa(S)=(1343)\)} \\
\hline
 \(11\) & \(39\) &  \(7\) &   \(2\) & \(1/25\) & \(117\) &    \(1,6,8\) & \(26\)       & \(3\)   &  \(3/27\) & \((4211)^3\) & \((3111)^3\) \\
 \(11\) &   \(\) &   \(\) &    \(\) &     \(\) & \(138\) &    \(2,6,7\) & \(23\)       & \(3\)   &  \(3/27\) & \((4211)^3\) & \((3111)^3\) \\
 \(11\) &   \(\) &   \(\) &    \(\) &     \(\) & \(162\) &    \(3,8;4\) & \(29;32\)    & \(3;4\) &  \(3/27\) & \((5211)^3\) & \((3111)^3\) \\
 \(11\) &   \(\) &   \(\) &    \(\) &     \(\) & \(184\) & \(21,23,25\) & \(23\)       & \(3\)   &  \(3/27\) & \((5211)^3\) & \((3111)^3\) \\

 \(11\) & \(39\) &  \(7\) &   \(9\) & \(1/25\) & \(117\) & \(21,23,25\) & \(26\)       & \(3\)   &  \(3/27\) & \((4211)^3\) & \((3111)^3\) \\
 \(11\) &   \(\) &   \(\) &    \(\) &     \(\) & \(138\) & \(10,14,18\) & \(23\)       & \(3\)   &  \(3/27\) & \((4211)^3\) & \((3111)^3\) \\
 \(11\) &   \(\) &   \(\) &    \(\) &     \(\) & \(162\) & \(11,16;15\) & \(29;32\)    & \(3;4\) &  \(3/27\) & \((5211)^3\) & \((3111)^3\) \\
 \(11\) &   \(\) &   \(\) &    \(\) &     \(\) & \(184\) &    \(1,6,8\) & \(23\)       & \(3\)   &  \(3/27\) & \((5211)^3\) & \((3111)^3\) \\

 \(11\) & \(44\) &  \(1\) &   \(3\) & \(1/25\) & \(114\) & \(11,15,16\) & \(26\)       & \(3\)   &  \(3/27\) & \((4211)^3\) & \((3111)^3\) \\
 \(11\) &   \(\) &   \(\) &    \(\) &     \(\) & \(136\) & \(21,23,25\) & \(23\)       & \(3\)   &  \(3/27\) & \((4211)^3\) & \((3111)^3\) \\
 \(11\) &   \(\) &   \(\) &    \(\) &     \(\) & \(159\) & \(21;23,25\) & \(29;32\)    & \(3;4\) &  \(3/27\) & \((5211)^3\) & \((3111)^3\) \\
 \(11\) &   \(\) &   \(\) &    \(\) &     \(\) & \(189\) & \(19,23,27\) & \(23\)       & \(3\)   &  \(3/27\) & \((5211)^3\) & \((3111)^3\) \\

 \(11\) & \(44\) &  \(1\) &   \(8\) & \(1/25\) & \(114\) &    \(3,4,8\) & \(26\)       & \(3\)   &  \(3/27\) & \((4211)^3\) & \((3111)^3\) \\
 \(11\) &   \(\) &   \(\) &    \(\) &     \(\) & \(136\) &    \(1,6,8\) & \(23\)       & \(3\)   &  \(3/27\) & \((4211)^3\) & \((3111)^3\) \\
 \(11\) &   \(\) &   \(\) &    \(\) &     \(\) & \(159\) &    \(1;6,8\) & \(29;32\)    & \(3;4\) &  \(3/27\) & \((5211)^3\) & \((3111)^3\) \\
 \(11\) &   \(\) &   \(\) &    \(\) &     \(\) & \(189\) &    \(3,4,8\) & \(23\)       & \(3\)   &  \(3/27\) & \((5211)^3\) & \((3111)^3\) \\
\hline
 \(11\) & \(54\) &  \(8\) &   \(2\) & \(1/25\) & \(120\) & \(13,17;12\) & \(26;29\)    & \(4\)   &  \(3/27\) & \((4211)^3\) & \((3111)^3\) \\
 \(11\) &   \(\) &   \(\) &    \(\) &     \(\) & \(140\) & \(20,22,27\) & \(29\)       & \(4\)   &  \(3/81\) & \((4211)^3\) & \((3111)^3\) \\
 \(11\) &   \(\) &   \(\) &    \(\) &     \(\) & \(140\) & \(38,40,45\) & \(29\)       & \(4\)   &  \(3/81\) & \((4211)^3\) & \((3111)^3\) \\
 \(11\) &   \(\) &   \(\) &    \(\) &     \(\) & \(140\) & \(56,58,63\) & \(29\)       & \(4\)   &  \(3/81\) & \((4211)^3\) & \((3111)^3\) \\
 \(11\) &   \(\) &   \(\) &    \(\) &     \(\) & \(166\) & \(20,22,27\) & \(23\)       & \(3\)   &  \(3/27\) & \((4211)^3\) & \((4111)^3\) \\
 \(11\) &   \(\) &   \(\) &    \(\) &     \(\) & \(176\) &    \(2,6,7\) & \(29\)       & \(4\)   &  \(3/27\) & \((4211)^3\) & \((4111)^3\) \\

 \(11\) & \(54\) &  \(8\) &   \(4\) & \(1/25\) & \(120\) & \(20,25;24\) & \(26;29\)    & \(4\)   &  \(3/27\) & \((4211)^3\) & \((3111)^3\) \\
 \(11\) &   \(\) &   \(\) &    \(\) &     \(\) & \(140\) & \(10,15,17\) & \(29\)       & \(4\)   &  \(3/81\) & \((4211)^3\) & \((3111)^3\) \\
 \(11\) &   \(\) &   \(\) &    \(\) &     \(\) & \(140\) & \(28,33,35\) & \(29\)       & \(4\)   &  \(3/81\) & \((4211)^3\) & \((3111)^3\) \\
 \(11\) &   \(\) &   \(\) &    \(\) &     \(\) & \(140\) & \(73,78,80\) & \(29\)       & \(4\)   &  \(3/81\) & \((4211)^3\) & \((3111)^3\) \\
 \(11\) &   \(\) &   \(\) &    \(\) &     \(\) & \(166\) & \(10,15,17\) & \(26\)       & \(3\)   &  \(3/27\) & \((4211)^3\) & \((4111)^3\) \\
 \(11\) &   \(\) &   \(\) &    \(\) &     \(\) & \(176\) & \(21,22,26\) & \(23\)       & \(3\)   &  \(3/27\) & \((4211)^3\) & \((4111)^3\) \\

 \(11\) & \(54\) &  \(8\) &   \(6\) & \(1/25\) & \(120\) & \(19,27;23\) & \(26;35\)    & \(4\)   &  \(3/27\) & \((4211)^3\) & \((3111)^3\) \\
 \(11\) &   \(\) &   \(\) &    \(\) &     \(\) & \(140\) & \(12,14,16\) & \(29\)       & \(4\)   &  \(3/81\) & \((4211)^3\) & \((3111)^3\) \\
 \(11\) &   \(\) &   \(\) &    \(\) &     \(\) & \(140\) & \(30,32,34\) & \(29\)       & \(4\)   &  \(3/81\) & \((4211)^3\) & \((3111)^3\) \\
 \(11\) &   \(\) &   \(\) &    \(\) &     \(\) & \(140\) & \(75,77,79\) & \(29\)       & \(4\)   &  \(3/81\) & \((4211)^3\) & \((3111)^3\) \\
 \(11\) &   \(\) &   \(\) &    \(\) &     \(\) & \(166\) & \(12,14,16\) & \(26\)       & \(3\)   &  \(3/27\) & \((4211)^3\) & \((4111)^3\) \\
 \(11\) &   \(\) &   \(\) &    \(\) &     \(\) & \(176\) & \(12,13,17\) & \(23\)       & \(3\)   &  \(3/27\) & \((4211)^3\) & \((4111)^3\) \\

 \(11\) & \(54\) &  \(8\) &   \(8\) & \(1/25\) & \(120\) &    \(3,4;8\) & \(26;32\)    & \(4\)   &  \(3/27\) & \((4211)^3\) & \((3111)^3\) \\
 \(11\) &   \(\) &   \(\) &    \(\) &     \(\) & \(140\) &    \(2,4,9\) & \(29\)       & \(4\)   &  \(3/81\) & \((4211)^3\) & \((3111)^3\) \\
 \(11\) &   \(\) &   \(\) &    \(\) &     \(\) & \(140\) & \(47,49,54\) & \(29\)       & \(4\)   &  \(3/81\) & \((4211)^3\) & \((3111)^3\) \\
 \(11\) &   \(\) &   \(\) &    \(\) &     \(\) & \(140\) & \(65,67,72\) & \(29\)       & \(4\)   &  \(3/81\) & \((4211)^3\) & \((3111)^3\) \\
 \(11\) &   \(\) &   \(\) &    \(\) &     \(\) & \(166\) &    \(2,4,9\) & \(23\)       & \(3\)   &  \(3/27\) & \((4211)^3\) & \((4111)^3\) \\
 \(11\) &   \(\) &   \(\) &    \(\) &     \(\) & \(176\) &    \(1,5,9\) & \(26\)       & \(3\)   &  \(3/27\) & \((4211)^3\) & \((4111)^3\) \\
\hline
\end{tabular}
\end{center}
\end{table}

%--------------------------------------------------------------------------------

\begin{theorem}
\label{thm:PeriF12Category3}
For an imaginary quadratic field \(K=\mathbb{Q}(\sqrt{d})\)
with discriminant \(d<0\),
bicyclic \(3\)-class group \(\mathrm{Cl}_3(K)\simeq C_3\times C_3\),
second \(3\)-class group \(M=\mathrm{G}_3^2(K)\), \(\mathrm{ord}(M)=3^{11}\),
capitulation type \(\mathrm{F}.12\), \(\varkappa(K)\sim (1343)\),
and any among the second abelian type invariants of third category,
\begin{equation}
\label{eqn:PeriF12Category3a}
\alpha^{(2)}(K)=\lbrack 11;(43;3321,(4211)^3),(32;3321,(3111)^3),(111;3321,(221)^3,(211)^9)^2\rbrack,
\end{equation}
\begin{equation}
\label{eqn:PeriF12Category3b}
\alpha^{(2)}(K)=\lbrack 11;(43;3321,(5211)^3),(32;3321,(3111)^3),(111;3321,(221)^3,(211)^9)^2\rbrack,
\end{equation}
\begin{equation}
\label{eqn:PeriF12Category3c}
\alpha^{(2)}(K)=\lbrack 11;(43;3321,(4211)^3),(32;3321,(4111)^3),(111;3321,(221)^3,(211)^9)^2\rbrack,
\end{equation}
the \(3\)-class field tower \(\mathrm{F}_3^\infty(K)\)
has two possible lengths \(\ell_3(K)\in\lbrace 3,4\rbrace\).
\end{theorem}

\begin{proof}
This is an immediate consequence of Table
\ref{tbl:CorrespondencePeriCC4F12},
since for each of three patterns in Formulas
\eqref{eqn:PeriF12Category3a},
\eqref{eqn:PeriF12Category3b},
\eqref{eqn:PeriF12Category3c}
there exist Schur \(\sigma\)-groups \(S\) with soluble lengths \(\mathrm{sl}(S)\in\lbrace 3,4\rbrace\).
Incidentally, 
\end{proof}

%--------------------------------------------------------------------------------

\begin{example}
\label{exm:PeriF12Category3}
According to Table
\ref{tbl:ATI2cc4Peri},
the quadratic fields \(K=\mathbb{Q}(\sqrt{d})\) with discriminants
\(d\in\lbrace -278\,427,-382\,123\rbrace\)
have the capitulation type \(\mathrm{F}.12\), \(\varkappa(K)=(1343)\),
and \(2^{\text{nd}}\) ATI of third category
\eqref{eqn:PeriF12Category3a} for \(d=-382\,123\),
\eqref{eqn:PeriF12Category3b} for \(d=-278\,427\).
According to Theorem
\ref{thm:PeriF12Category3},
no sharp decision about the length of the \(3\)-class field tower \(\mathrm{F}_3^\infty(K)\)
can be drawn, since there are two possibilities \(\ell_3(K)\in\lbrace 3,4\rbrace\).
However, the skeleton type is
\(\mathrm{d}.19\), \(j_1\in\lbrace 39,44\rbrace\), for
\eqref{eqn:PeriF12Category3b}.
%and \(\mathrm{d}.23\), \(j_1=54\), for
%\eqref{eqn:PeriF13Category3c}.
\end{example}

%\newpage
%--------------------------------------------------------------------------------

\renewcommand{\arraystretch}{1.0}

\begin{table}[hb]
\caption{Correspondence for periodic \(M=S/S^{\prime\prime}\), \(\mathrm{cc}(M)=4\), of type \(\mathrm{F}.11\)}
\label{tbl:CorrespondencePeriCC4F11}
\begin{center}
\begin{tabular}{|c|c|c|c|c||c|c|c|c|c||c|c|}
\hline
   Type &  &  &  &  &  &  &  &  &  & \multicolumn{2}{|c|}{\(\alpha^{(2)}(S)\)} \\
 lo\((M)\)& \(j_1\) & \(j_2\) & \(j_3\) & \(\#\) & \(k_1\) & \(k_3\) & lo\((S)\) & sl\((S)\) & \(\#\) & \(A_1\) & \(A_2\) \\
\hline
 \(\mathrm{F}.11\) &  &  &  &  &  &  & \multicolumn{4}{|c|}{\(\varkappa(S)=(1143)\)} \\
\hline
 \(11\) & \(54\) &  \(8\) &   \(1\) & \(1/25\) & \(120\) & \(15,16;11\) & \(26;29\)    & \(4\)   &  \(3/27\) & \((5211)^3\) & \((3111)^3\) \\
 \(11\) &   \(\) &   \(\) &    \(\) &     \(\) & \(140\) & \(19,24,26\) & \(29\)       & \(3\)   &  \(3/81\) & \((5211)^3\) & \((3111)^3\) \\
 \(11\) &   \(\) &   \(\) &    \(\) &     \(\) & \(140\) & \(37,42,44\) & \(29\)       & \(3\)   &  \(3/81\) & \((5211)^3\) & \((3111)^3\) \\
 \(11\) &   \(\) &   \(\) &    \(\) &     \(\) & \(140\) & \(55,60,62\) & \(29\)       & \(3\)   &  \(3/81\) & \((5211)^3\) & \((3111)^3\) \\
 \(11\) &   \(\) &   \(\) &    \(\) &     \(\) & \(166\) & \(19,24,26\) & \(23\)       & \(3\)   &  \(3/27\) & \((4211)^3\) & \((4111)^3\) \\
 \(11\) &   \(\) &   \(\) &    \(\) &     \(\) & \(176\) & \(20,24,25\) & \(23\)       & \(3\)   &  \(3/27\) & \((4211)^3\) & \((4111)^3\) \\

 \(11\) & \(54\) &  \(8\) &   \(9\) & \(1/25\) & \(120\) &    \(1,5;9\) & \(26;32\)    & \(4\)   &  \(3/27\) & \((5211)^3\) & \((3111)^3\) \\
 \(11\) &   \(\) &   \(\) &    \(\) &     \(\) & \(140\) &    \(3,5,7\) & \(29\)       & \(3\)   &  \(3/81\) & \((5211)^3\) & \((3111)^3\) \\
 \(11\) &   \(\) &   \(\) &    \(\) &     \(\) & \(140\) & \(48,50,52\) & \(29\)       & \(3\)   &  \(3/81\) & \((5211)^3\) & \((3111)^3\) \\
 \(11\) &   \(\) &   \(\) &    \(\) &     \(\) & \(140\) & \(66,68,70\) & \(29\)       & \(3\)   &  \(3/81\) & \((5211)^3\) & \((3111)^3\) \\
 \(11\) &   \(\) &   \(\) &    \(\) &     \(\) & \(166\) &    \(3,5,7\) & \(23\)       & \(3\)   &  \(3/27\) & \((4211)^3\) & \((4111)^3\) \\
 \(11\) &   \(\) &   \(\) &    \(\) &     \(\) & \(176\) & \(10,14,18\) & \(23\)       & \(3\)   &  \(3/27\) & \((4211)^3\) & \((4111)^3\) \\
\hline
 \(11\) & \(57\) &  \(1\) &   \(1\) & \(1/15\) & \(124\) &    \(1,6,8\) & \(23\)       & \(3\)   &  \(3/27\) & \((4211)^3\) & \((4111)^3\) \\
 \(11\) &   \(\) &   \(\) &    \(\) &     \(\) & \(124\) & \(21,23,25\) & \(23\)       & \(3\)   &  \(3/27\) & \((4211)^3\) & \((4111)^3\) \\
 \(11\) &   \(\) &   \(\) &    \(\) &     \(\) & \(168\) &   \(1,14;8\) & \(26;29\)    & \(4\)   &  \(3/15\) & \((5211)^3\) & \((3111)^3\) \\
 \(11\) &   \(\) &   \(\) &    \(\) &     \(\) & \(197\) &   \(3,4,10\) & \(29\)       & \(3\)   &  \(3/45\) & \((5211)^3\) & \((3111)^3\) \\
 \(11\) &   \(\) &   \(\) &    \(\) &     \(\) & \(197\) & \(14,21,30\) & \(29\)       & \(3\)   &  \(3/45\) & \((5211)^3\) & \((3111)^3\) \\
 \(11\) &   \(\) &   \(\) &    \(\) &     \(\) & \(197\) & \(34,36,43\) & \(29\)       & \(3\)   &  \(3/45\) & \((5211)^3\) & \((3111)^3\) \\ 
 
 \(11\) & \(59\) &  \(6\) &   \(1\) & \(1/15\) & \(129\) & \(10,14,18\) & \(23\)       & \(3\)   &  \(3/27\) & \((4211)^3\) & \((4111)^3\) \\
 \(11\) &   \(\) &   \(\) &    \(\) &     \(\) & \(129\) & \(21,22,26\) & \(23\)       & \(3\)   &  \(3/27\) & \((4211)^3\) & \((4111)^3\) \\
 \(11\) &   \(\) &   \(\) &    \(\) &     \(\) & \(173\) &  \(1,10;13\) & \(26;32\)    & \(4\)   &  \(3/27\) & \((5211)^3\) & \((3111)^3\) \\
 \(11\) &   \(\) &   \(\) &    \(\) &     \(\) & \(198\) &   \(2,4,12\) & \(29\)       & \(3\)   &  \(3/45\) & \((5211)^3\) & \((3111)^3\) \\
 \(11\) &   \(\) &   \(\) &    \(\) &     \(\) & \(198\) & \(18,21,29\) & \(29\)       & \(3\)   &  \(3/45\) & \((5211)^3\) & \((3111)^3\) \\ 
 \(11\) &   \(\) &   \(\) &    \(\) &     \(\) & \(198\) & \(33,39,44\) & \(29\)       & \(3\)   &  \(3/45\) & \((5211)^3\) & \((3111)^3\) \\
\hline
\end{tabular}
\end{center}
\end{table}

%--------------------------------------------------------------------------------

\begin{theorem}
\label{thm:PeriF11Category3}
An imaginary quadratic field \(K=\mathbb{Q}(\sqrt{d})\)
with fundamental discriminant \(d<0\),
elementary bicyclic \(3\)-class group \(\mathrm{Cl}_3(K)\),
second \(3\)-class group \(M=\mathrm{G}_3^2(K)\), \(\mathrm{ord}(M)=3^{11}\),
capitulation type \(\mathrm{F}.11\), \(\varkappa(K)\sim (1143)\),
and second abelian type invariants of third category
\begin{equation}
\label{eqn:PeriF11Category3a}
\alpha^{(2)}(K)=\lbrack 11;(43;3321,(4211)^3),(32;3321,(4111)^3),(111;3321,(221)^3,(211)^9)^2\rbrack
\end{equation}
has a \(3\)-class field tower \(\mathrm{F}_3^\infty(K)\) of precise length \(\ell_3(K)=3\).
For second abelian type invariants
\begin{equation}
\label{eqn:PeriF11Category3b}
\alpha^{(2)}(K)=\lbrack 11;(43;3321,(5211)^3),(32;3321,(3111)^3),(111;3321,(221)^3,(211)^9)^2\rbrack,
\end{equation}
the \(3\)-class field tower \(\mathrm{F}_3^\infty(K)\)
has two possible lengths \(\ell_3(K)\in\lbrace 3,4\rbrace\).
\end{theorem}

\begin{proof}
This is an immediate consequence of Table
\ref{tbl:CorrespondencePeriCC4F11},
since the settled roots \(R\) arising from the
unique non-metabelian \(B-\#4;k_1\) with \(k_1\in\lbrace 120,168,173\rbrace\)
which lead to Schur \(\sigma\)-groups \(S\) with soluble length \(\mathrm{sl}(S)=4\)
have the second abelian type invariants in Formula
\eqref{eqn:PeriF11Category3b}.
\end{proof}

%--------------------------------------------------------------------------------

\begin{example}
\label{exm:PeriF11Category3}
According to Table
\ref{tbl:ATI2cc4Peri},
the quadratic field \(K=\mathbb{Q}(\sqrt{d})\) with discriminant
\(d=-469\,787\)
has the capitulation type \(\mathrm{F}.11\), \(\varkappa(K)=(1143)\),
and second abelian type invariants of third category
\eqref{eqn:PeriF11Category3b}.
According to Theorem
\ref{thm:PeriF11Category3},
no sharp decision about the length of the \(3\)-class field tower \(\mathrm{F}_3^\infty(K)\)
can be drawn, since there are two possibilities \(\ell_3(K)\in\lbrace 3,4\rbrace\).
However, outside of the range in Table
\ref{tbl:ATI2cc4Peri},
we found quadratic fields \(K=\mathbb{Q}(\sqrt{d})\) with discriminants
\(d\in\lbrace -2\,005\,880,-2\,554\,868,\)
\(-3\,283\,223,-3\,937\,999\rbrace\),
capitulation type \(\mathrm{F}.11\), \(\varkappa(K)=(1143)\),
and second abelian type invariants of third category
\eqref{eqn:PeriF11Category3a}.
According to Theorem
\ref{thm:PeriF11Category3},
they have a \(3\)-class field tower \(\mathrm{F}_3^\infty(K)\) of precise length \(\ell_3(K)=3\).
\end{example}

%--------------------------------------------------------------------------------

\begin{theorem}
\label{thm:PeriF7Category3}
For an imaginary quadratic field \(K=\mathbb{Q}(\sqrt{d})\)
with discriminant \(d<0\),
elementary bicyclic \(3\)-class group \(\mathrm{Cl}_3(K)\),
second \(3\)-class group \(M=\mathrm{G}_3^2(K)\), \(\mathrm{ord}(M)=3^{11}\),
and capitulation type \(\mathrm{F}.7\), \(\varkappa(K)\sim (3443)\),
the second abelian type invariants of third category
must have the shape
\begin{equation}
\label{eqn:PeriF7Category3}
\alpha^{(2)}(K)=\lbrack 11;(43;3321,(4111)^3),(32;3321,(3111)^3),(111;3321,(221)^3,(211)^9)^2\rbrack,
\end{equation}
and the \(3\)-class field tower \(\mathrm{F}_3^\infty(K)\)
has two possible lengths \(\ell_3(K)\in\lbrace 3,4\rbrace\).
\end{theorem}

\begin{proof}
This is an immediate consequence of Table
\ref{tbl:CorrespondencePeriCC4F7},
since the settled roots \(R\) arising from the
unique non-metabelian \(B-\#4;k_1\) with \(k_1\in\lbrace 136,138,184\rbrace\)
which lead to Schur \(\sigma\)-groups \(S\) with soluble length \(\mathrm{sl}(S)=4\)
also have the common second abelian type invariants in Formula
\eqref{eqn:PeriF7Category3}.
\end{proof}

%\newpage
%--------------------------------------------------------------------------------

\renewcommand{\arraystretch}{1.0}

\begin{table}[hb]
\caption{Correspondence for periodic \(M=S/S^{\prime\prime}\), \(\mathrm{cc}(M)=4\), of type \(\mathrm{F}.7\)}
\label{tbl:CorrespondencePeriCC4F7}
\begin{center}
\begin{tabular}{|c|c|c|c|c||c|c|c|c|c||c|c|}
\hline
   Type &  &  &  &  &  &  &  &  &  & \multicolumn{2}{|c|}{\(\alpha^{(2)}(S)\)} \\
 lo\((M)\)& \(j_1\) & \(j_2\) & \(j_3\) & \(\#\) & \(k_1\) & \(k_3\) & lo\((S)\) & sl\((S)\) & \(\#\) & \(A_1\) & \(A_2\) \\
\hline
 \(\mathrm{F}.7\) &  &  &  &  &  &  &  &  &  & \multicolumn{2}{|c|}{\(\varkappa(S)=(3443)\)} \\
\hline
 \(11\) & \(39\) &  \(7\) &   \(5\) & \(1/25\) & \(117\) &    \(2,4,9\) & \(26\)       & \(3\)   &  \(3/27\) & \((4211)^3\) & \((3111)^3\) \\
 \(11\) &   \(\) &   \(\) &    \(\) &     \(\) & \(138\) &    \(3,4,8\) & \(26\)       & \(3\)   &  \(3/27\) & \((4211)^3\) & \((3111)^3\) \\
 \(11\) &   \(\) &   \(\) &    \(\) &     \(\) & \(162\) &    \(1,5,9\) & \(26\)       & \(3\)   &  \(3/27\) & \((4211)^3\) & \((3111)^3\) \\
 \(11\) &   \(\) &   \(\) &    \(\) &     \(\) & \(184\) & \(12,14,16\) & \(26\)       & \(3\)   &  \(3/27\) & \((4211)^3\) & \((3111)^3\) \\

 \(11\) & \(39\) &  \(7\) &   \(6\) & \(1/25\) & \(117\) & \(20,22,27\) & \(26\)       & \(3\)   &  \(3/27\) & \((4211)^3\) & \((3111)^3\) \\
 \(11\) &   \(\) &   \(\) &    \(\) &     \(\) & \(138\) & \(12,17;13\) & \(26;29\)    & \(3;4\) &  \(3/27\) & \((4211)^3\) & \((3111)^3\) \\
 \(11\) &   \(\) &   \(\) &    \(\) &     \(\) & \(162\) & \(10,14,18\) & \(26\)       & \(3\)   &  \(3/27\) & \((4211)^3\) & \((3111)^3\) \\
 \(11\) &   \(\) &   \(\) &    \(\) &     \(\) & \(184\) & \(10,15,17\) & \(29\)       & \(4\)   &  \(3/27\) & \((4211)^3\) & \((3111)^3\) \\

 \(11\) & \(44\) &  \(1\) &   \(5\) & \(1/25\) & \(114\) &    \(1,5,9\) & \(26\)       & \(3\)   &  \(3/27\) & \((4211)^3\) & \((3111)^3\) \\
 \(11\) &   \(\) &   \(\) &    \(\) &     \(\) & \(136\) &    \(2;4,9\) & \(26;29,32\) & \(3;4\) &  \(3/27\) & \((4211)^3\) & \((3111)^3\) \\
 \(11\) &   \(\) &   \(\) &    \(\) &     \(\) & \(159\) &    \(2,4,9\) & \(26\)       & \(3\)   &  \(3/27\) & \((4211)^3\) & \((3111)^3\) \\
 \(11\) &   \(\) &   \(\) &    \(\) &     \(\) & \(189\) & \(12,13,17\) & \(26\)       & \(3\)   &  \(3/27\) & \((4211)^3\) & \((3111)^3\) \\

 \(11\) & \(44\) &  \(1\) &   \(6\) & \(1/25\) & \(114\) & \(10,14,18\) & \(26\)       & \(3\)   &  \(3/27\) & \((4211)^3\) & \((3111)^3\) \\
 \(11\) &   \(\) &   \(\) &    \(\) &     \(\) & \(136\) & \(20,27;22\) & \(26;29\)    & \(3;4\) &  \(3/27\) & \((4211)^3\) & \((3111)^3\) \\
 \(11\) &   \(\) &   \(\) &    \(\) &     \(\) & \(159\) & \(20,22,27\) & \(26\)       & \(3\)   &  \(3/27\) & \((4211)^3\) & \((3111)^3\) \\
 \(11\) &   \(\) &   \(\) &    \(\) &     \(\) & \(189\) & \(10,14,18\) & \(26\)       & \(3\)   &  \(3/27\) & \((4211)^3\) & \((3111)^3\) \\
\hline
\end{tabular}
\end{center}
\end{table}

%--------------------------------------------------------------------------------

\begin{example}
\label{exm:PeriF7Category3}
According to Table
\ref{tbl:ATI2cc4Peri},
the quadratic fields \(K=\mathbb{Q}(\sqrt{d})\) with discriminants
\(d\in\lbrace -643\,011,-797\,556\rbrace\)
have the capitulation type \(\mathrm{F}.7\), \(\varkappa(K)=(3443)\),
and second abelian type invariants of third category
\eqref{eqn:PeriF7Category3}.
According to Theorem
\ref{thm:PeriF7Category3},
no sharp decision about the length of the \(3\)-class field tower \(\mathrm{F}_3^\infty(K)\)
can be drawn, since there are two possibilities \(\ell_3(K)\in\lbrace 3,4\rbrace\).
\end{example}

%--------------------------------------------------------------------------------

\begin{corollary}
\label{cor:PeriodicCc4Category3}
The logarithmic order
of the Schur \(\sigma\)-group \(S=\mathrm{Gal}(\mathrm{F}_3^\infty/K)\)
is \(\mathrm{lo}(S)\in\lbrace 23,26,29\rbrace\) for a \(3\)-class field tower of length \(\ell_3(K)=3\),
and it depends on the capitulation type \(\varkappa(K)\),
\[
\mathrm{lo}(S)\in
\begin{cases}
\lbrace 26,29,32,35,41\rbrace & \text{ for type } \mathrm{F}.13, \\
\lbrace 26,29,32,35\rbrace & \text{ for type } \mathrm{F}.12, \\
\lbrace 26,29,32\rbrace & \text{ for type } \mathrm{F}.11, \\
\lbrace 29,32\rbrace & \text{ for type } \mathrm{F}.7,
\end{cases}
\] 
in the case of a tower of length \(\ell_3(K)=4\).
\end{corollary}

\begin{proof}
This is a common consequence of the Tables
\ref{tbl:CorrespondencePeriCC4F13},
\ref{tbl:CorrespondencePeriCC4F12},
\ref{tbl:CorrespondencePeriCC4F11},
\ref{tbl:CorrespondencePeriCC4F7}.
\end{proof}

%\newpage
%--------------------------------------------------------------------------------

\section{Realization of Sporadic Metabelian \(\sigma\)-Groups \(\mathrm{G}_3^2(K)\) of Coclass \(6\)}
\label{s:SporCc6}

\noindent
In Table
\ref{tbl:ATI2cc6Spor},
second abelian type invariants \(\alpha^{(2)}(K)\) of imaginary quadratic fields \(K=\mathbb{Q}(\sqrt{d})\)
with discriminants in the most extensive range \(-10^7<d<0\),
bicyclic \(3\)-class group \(\mathrm{Cl}_3(K)\simeq C_3\times C_3\),
second \(3\)-class group \(M=\mathrm{G}_3^2(K)\), \(\mathrm{ord}(M)=3^{13}\),
and type \(\mathrm{F}\) are given,
ordered by capitulation subtypes \(\varkappa(K)\).
\textbf{Boldface} font indicates the categories \(1\) and \(2\).

%\newpage
%--------------------------------------------------------------------------------

\renewcommand{\arraystretch}{1.0}

\begin{table}[ht]
\caption{\(2^{\text{nd}}\) ATI of \(K=\mathbb{Q}(\sqrt{d})\), \(d<0\), with sporadic \(M=\mathrm{G}_3^2(K)\), \(\mathrm{cc}(M)=6\)}
\label{tbl:ATI2cc6Spor}
\begin{center}
\begin{tabular}{|r|c|c|c|c|}
\hline
         Type & \multicolumn{4}{|c|}{\(\alpha^{(2)}(K)=\lbrack 11;(43;3332,A_1),(43;3332,A_2),(111;3332,A_3),(111;3332,A_4)\rbrack\)} \\
       \(-d\) & \(A_1\) & \(A_2\) & \(A_3\) & \(A_4\) \\
\hline
 \(\mathrm{F}.7\) & \multicolumn{4}{|c|}{\(\varkappa(K)=(3443)\)} \\
\hline
 \(1\,677\,768\) & \((5211)^3\)           & \((4211)^3\)           & \(\mathbf{(3211)^3,(2111)^6},(221)^3\)    & \(\mathbf{(2111)^6,(1111)^6}\)          \\
 \(5\,053\,191\) & \((4211)^3\)           & \((4211)^3\)           & \((221)^3,(211)^9\)                       & \((221)^3,(211)^9\)                     \\
 \(8\,723\,023\) & \((4211)^3\)           & \((4211)^3\)           & \((221)^3,(211)^9\)                       & \((221)^3,(211)^9\)                     \\
\hline
 \(\mathrm{F}.11\) & \multicolumn{4}{|c|}{\(\varkappa(K)=(1143)\)} \\
\hline
 \(4\,838\,891\) & \((5211)^3\)           & \((4211)^3\)           & \((221)^3,(211)^9\)                       & \((221)^3,(211)^9\)                     \\
 \(5\,427\,023\) & \((4211)^3\)           & \((4211)^3\)           & \(\mathbf{(2111)^9},(221)^3\)             & \(\mathbf{(1^5)^3,(2111)^3,(1111)^6}\). \\
 \(8\,493\,815\) & \((5211)^3\)           & \((4211)^3\)           & \((221)^3,(211)^9\)                       & \((221)^3,(211)^9\)                     \\
\hline
 \(\mathrm{F}.12\) & \multicolumn{4}{|c|}{\(\varkappa(K)=(1343)\)} \\
\hline
    \(423\,640\) & \((4211)^3\)           & \((4211)^3\)           & \((221)^3,(211)^9\)                       & \((221)^3,(211)^9\)                     \\
 \(8\,751\,215\) & \((5211)^3\)           & \((4211)^3\)           & \((221)^3,(211)^9\)                       & \((221)^3,(211)^9\)                     \\
\hline
 \(\mathrm{F}.13\) & \multicolumn{4}{|c|}{\(\varkappa(K)=(3143)\)} \\
\hline
 \(2\,383\,059\) & \(\mathbf{(4221)^3}\)  & \(\mathbf{(4221)^3}\)  & \(\mathbf{(21^4)^3,(1^5)^6,(2211)^3}\)    & \(\mathbf{(21^4)^3,(2221)^6,(2211)^3}\) \\
 \(5\,444\,651\) & \((4211)^3\)           & \((4211)^3\)           & \(\mathbf{(2211)^3,(2111)^6},(221)^3\)    & \(\mathbf{(2211)^3,(2111)^3,(1111)^6}\) \\
 \(5\,606\,283\) & \((4211)^3\)           & \((4211)^3\)           & \(\mathbf{(2111)^9},(221)^3\)             & \(\mathbf{(2111)^6,(1111)^6}\)          \\
 \(5\,765\,812\) & \(\mathbf{(5221)^3}\)  & \(\mathbf{(4221)^3}\)  & \(\mathbf{(21^4)^3,(1^5)^6,(2211)^3}\)    & \(\mathbf{(21^4)^3,(2211)^9}\)          \\
 \(6\,863\,219\) & \((5211)^3\)           & \((4211)^3\)           & \((221)^3,(211)^9\)                       & \((221)^3,(211)^9\)                     \\
 \(8\,963\,839\) & \((4211)^3\)           & \((4211)^3\)           & \(\mathbf{(2211)^3,(2111)^6},(221)^3\)    & \(\mathbf{(1^5)^3,(2111)^3,(1111)^6}\)  \\
\hline
\end{tabular}
\end{center}
\end{table}

%\newpage
%--------------------------------------------------------------------------------

\bigskip
\noindent
Throughout the sequel, we restrict our investigations to the tame situation of category \(3\).

%--------------------------------------------------------------------------------

\bigskip
\noindent
In the Tables
\ref{tbl:CorrespondenceSporCC6F7},
\ref{tbl:CorrespondenceSporCC6F13},
\ref{tbl:CorrespondenceSporCC6F12},
\ref{tbl:CorrespondenceSporCC6F11},
the Schur \(\sigma\)-groups \(S\)
have the settled root \(R=B-\#4;k_1-\#2;k_2-\#4;k_3\)
with some unique relative identifier \(k_2=k_2(k_1)\),
metabelianization \(M\simeq B-\#2;33-\#2;25-\#2;j\),
and \(2^{\text{nd}}\) AQI 
\(\alpha^{(2)}(S)=\lbrack 11;(43;3332,A_1),(43;3332,A_2),(111;3332,A_3),(111;3332,A_4)\rbrack\).

%--------------------------------------------------------------------------------

\renewcommand{\arraystretch}{1.0}

\begin{table}[ht]
\caption{Correspondence for sporadic \(M=S/S^{\prime\prime}\), \(\mathrm{cc}(M)=6\), of type \(\mathrm{F}.7\)}
\label{tbl:CorrespondenceSporCC6F7}
\begin{center}
\begin{tabular}{|c|c||c|c|c|c|c||c|c|c|c|}
\hline
   Type &  &  &  &  &  &  & \multicolumn{4}{|c|}{\(\alpha^{(2)}(S)\)} \\
 lo\((M)\)& \(j\) & \(k_1\) & \(k_3\) & lo\((S)\) & sl\((S)\) & \(\#\) & \(A_1\) & \(A_2\) & \(A_3\) & \(A_4\) \\
\hline
 \(\mathrm{F}.7\) &  &  &  &  &  &  & \multicolumn{4}{|c|}{\(\varkappa(S)=(3443)\)} \\
\hline
 \(13\) & \(59\) & \(148\) &  \(8\) & \(26\) & \(3\) & \(1/45\) & \((4211)^3\) & \((4211)^3\) & \((221)^3,(211)^9\) & \((221)^3,(211)^9\) \\
   \(\) &   \(\) &    \(\) & \(21\) & \(26\) & \(3\) & \(1/45\) & \((4211)^3\) & \((4211)^3\) & \((221)^3,(211)^9\) & \((221)^3,(211)^9\) \\
   \(\) &   \(\) & \(179\) &  \(8\) & \(26\) & \(3\) & \(1/45\) & \((4211)^3\) & \((4211)^3\) & \((221)^3,(211)^9\) & \((221)^3,(211)^9\) \\
   \(\) &   \(\) &    \(\) & \(42\) & \(26\) & \(3\) & \(1/45\) & \((4211)^3\) & \((4211)^3\) & \((221)^3,(211)^9\) & \((221)^3,(211)^9\) \\
 \(13\) & \(60\) & \(148\) & \(24\) & \(26\) & \(3\) & \(1/45\) & \((4211)^3\) & \((4211)^3\) & \((221)^3,(211)^9\) & \((221)^3,(211)^9\) \\
   \(\) &   \(\) &    \(\) & \(45\) & \(26\) & \(3\) & \(1/45\) & \((4211)^3\) & \((4211)^3\) & \((221)^3,(211)^9\) & \((221)^3,(211)^9\) \\
   \(\) &   \(\) & \(179\) & \(24\) & \(26\) & \(3\) & \(1/45\) & \((4211)^3\) & \((4211)^3\) & \((221)^3,(211)^9\) & \((221)^3,(211)^9\) \\
 \(13\) & \(62\) & \(148\) &  \(6\) & \(26\) & \(3\) & \(1/45\) & \((4211)^3\) & \((4211)^3\) & \((221)^3,(211)^9\) & \((221)^3,(211)^9\) \\
   \(\) &   \(\) &    \(\) & \(41\) & \(26\) & \(3\) & \(1/45\) & \((4211)^3\) & \((4211)^3\) & \((221)^3,(211)^9\) & \((221)^3,(211)^9\) \\
   \(\) &   \(\) & \(179\) &  \(6\) & \(26\) & \(3\) & \(1/45\) & \((4211)^3\) & \((4211)^3\) & \((221)^3,(211)^9\) & \((221)^3,(211)^9\) \\
\hline
\end{tabular}
\end{center}
\end{table}

%--------------------------------------------------------------------------------

\begin{theorem}
\label{thm:F7Category3Cc6}
Let \(K=\mathbb{Q}(\sqrt{d})\) be an imaginary quadratic field
with fundamental discriminant \(d<0\),
elementary bicyclic \(3\)-class group \(\mathrm{Cl}_3(K)\simeq C_3\times C_3\),
capitulation type \(\mathrm{F}.7\), \(\varkappa(K)\sim (3443)\),
and first abelian type invariants
\(\alpha^{(1)}(K)=\lbrack 11;43,43,111,111\rbrack\).
If the second abelian type invariants are of third category,
they must have the shape
\begin{equation}
\label{eqn:F7Category3Cc6}
\alpha^{(2)}(K)=\lbrack 11;(43;3332,(4211)^3),(43;3332,(4211)^3),(111;3332,(221)^3,(211)^9)^2\rbrack
\end{equation}
and the \(3\)-class field tower \(\mathrm{F}_3^\infty(K)\) has precise length \(\ell_3(K)=3\) and relative degree \(3^{26}\) over \(K\).
\end{theorem}

\begin{proof}
This is an immediate consequence of Table
\ref{tbl:CorrespondenceSporCC6F7},
since all the non-metabelian roots \(R=B-\#4;k_1\)
with \(k_1\in\lbrace 148,179\rbrace\)
lead to Schur \(\sigma\)-groups \(S\)
with transfer kernel type \(\mathrm{F}.7\), \(\varkappa(K)=(3443)\),
whose soluble length is \(\mathrm{sl}(S)=3\),
logarithmic order \(\mathrm{lo}(S)=26\),
and which have the second abelian type invariants in Formula
\eqref{eqn:F7Category3Cc6}.
\end{proof}

%--------------------------------------------------------------------------------

\begin{example}
\label{exm:F7Category3Cc6}
According to Table
\ref{tbl:ATI2cc6Spor},
the quadratic fields \(K=\mathbb{Q}(\sqrt{d})\)
with discriminants \(d\in\lbrace -5\,053\,191,-8\,723\,023\rbrace\)
satisfy the assumptions of Theorem
\ref{thm:F7Category3Cc6}.
Thus, their \(3\)-class field tower has precise length \(\ell_3(K)=3\)
and relative degree \(\lbrack\mathrm{F}_3^\infty(K):K\rbrack=3^{26}\).
\end{example}

%--------------------------------------------------------------------------------

\renewcommand{\arraystretch}{1.0}

\begin{table}[hb]
\caption{Correspondence for sporadic \(M=S/S^{\prime\prime}\), \(\mathrm{cc}(M)=6\), of type \(\mathrm{F}.13\)}
\label{tbl:CorrespondenceSporCC6F13}
\begin{center}
\begin{tabular}{|c|c||c|c|c|c|c||c|c|c|c|}
\hline
   Type &  &  &  &  &  &  & \multicolumn{4}{|c|}{\(\alpha^{(2)}(S)\)} \\
 lo\((M)\)& \(j\) & \(k_1\) & \(k_3\) & lo\((S)\) & sl\((S)\) & \(\#\) & \(A_1\) & \(A_2\) & \(A_3\) & \(A_4\) \\
\hline
 \(\mathrm{F}.13\) &  &  &  &  &  &  & \multicolumn{4}{|c|}{\(\varkappa(S)=(3143)\)} \\
\hline
 \(13\) & \(45\) & \(148\) &  \(7\) & \(26\) & \(3\) & \(1/45\) & \((5211)^3\) & \((4211)^3\) & \((221)^3,(211)^9\) & \((221)^3,(211)^9\) \\
   \(\) &   \(\) &    \(\) & \(37\) & \(26\) & \(3\) & \(1/45\) & \((5211)^3\) & \((4211)^3\) & \((221)^3,(211)^9\) & \((221)^3,(211)^9\) \\
   \(\) &   \(\) & \(179\) &  \(7\) & \(32\) & \(4\) & \(1/45\) & \((4211)^3\) & \((4211)^3\) & \((221)^3,(211)^9\) & \((221)^3,(211)^9\) \\
   \(\) &   \(\) &    \(\) & \(37\) & \(29\) & \(3\) & \(1/45\) & \((4211)^3\) & \((4211)^3\) & \((221)^3,(211)^9\) & \((221)^3,(211)^9\) \\
 \(13\) & \(51\) & \(148\) & \(20\) & \(26\) & \(3\) & \(1/45\) & \((5211)^3\) & \((4211)^3\) & \((221)^3,(211)^9\) & \((221)^3,(211)^9\) \\
   \(\) &   \(\) &    \(\) & \(30\) & \(26\) & \(3\) & \(1/45\) & \((5211)^3\) & \((4211)^3\) & \((221)^3,(211)^9\) & \((221)^3,(211)^9\) \\
   \(\) &   \(\) & \(179\) & \(20\) & \(29\) & \(3\) & \(1/45\) & \((4211)^3\) & \((4211)^3\) & \((221)^3,(211)^9\) & \((221)^3,(211)^9\) \\
   \(\) &   \(\) &    \(\) & \(30\) & \(32\) & \(4\) & \(1/45\) & \((4211)^3\) & \((4211)^3\) & \((221)^3,(211)^9\) & \((221)^3,(211)^9\) \\
 \(13\) & \(54\) & \(148\) &  \(3\) & \(26\) & \(3\) & \(1/45\) & \((5211)^3\) & \((4211)^3\) & \((221)^3,(211)^9\) & \((221)^3,(211)^9\) \\
   \(\) &   \(\) &    \(\) & \(40\) & \(26\) & \(3\) & \(1/45\) & \((5211)^3\) & \((4211)^3\) & \((221)^3,(211)^9\) & \((221)^3,(211)^9\) \\
   \(\) &   \(\) & \(179\) &  \(3\) & \(29\) & \(3\) & \(1/45\) & \((4211)^3\) & \((4211)^3\) & \((221)^3,(211)^9\) & \((221)^3,(211)^9\) \\
   \(\) &   \(\) &    \(\) & \(14\) & \(29\) & \(3\) & \(1/45\) & \((4211)^3\) & \((4211)^3\) & \((221)^3,(211)^9\) & \((221)^3,(211)^9\) \\
 \(13\) & \(56\) & \(148\) & \(17\) & \(26\) & \(3\) & \(1/45\) & \((5211)^3\) & \((4211)^3\) & \((221)^3,(211)^9\) & \((221)^3,(211)^9\) \\
   \(\) &   \(\) &    \(\) & \(44\) & \(26\) & \(3\) & \(1/45\) & \((5211)^3\) & \((4211)^3\) & \((221)^3,(211)^9\) & \((221)^3,(211)^9\) \\
   \(\) &   \(\) & \(179\) & \(17\) & \(29\) & \(3\) & \(1/45\) & \((4211)^3\) & \((4211)^3\) & \((221)^3,(211)^9\) & \((221)^3,(211)^9\) \\
   \(\) &   \(\) &    \(\) & \(19\) & \(29\) & \(3\) & \(1/45\) & \((4211)^3\) & \((4211)^3\) & \((221)^3,(211)^9\) & \((221)^3,(211)^9\) \\
\hline
\end{tabular}
\end{center}
\end{table}

%--------------------------------------------------------------------------------

\begin{theorem}
\label{thm:F13Category3Cc6}
Let \(K=\mathbb{Q}(\sqrt{d})\) be an imaginary quadratic field
with fundamental discriminant \(d<0\),
elementary bicyclic \(3\)-class group \(\mathrm{Cl}_3(K)\simeq C_3\times C_3\),
capitulation type \(\mathrm{F}.13\), \(\varkappa(K)\sim (3143)\),
and first abelian type invariants
\(\alpha^{(1)}(K)=\lbrack 11;43,43,111,111\rbrack\).
If the second abelian type invariants are of third category,
they must either have the shape
\begin{equation}
\label{eqn:F13Category3Cc6a}
\alpha^{(2)}(K)=\lbrack 11;(43;3332,(5211)^3),(43;3332,(4211)^3),(111;3332,(221)^3,(211)^9)^2\rbrack
\end{equation}
and the \(3\)-class field tower \(\mathrm{F}_3^\infty(K)\) has precise length \(\ell_3(K)=3\) and relative degree \(3^{26}\) over \(K\),
or the shape
\begin{equation}
\label{eqn:F13Category3Cc6b}
\alpha^{(2)}(K)=\lbrack 11;(43;3332,(4211)^3),(43;3332,(4211)^3),(111;3332,(221)^3,(211)^9)^2\rbrack
\end{equation}
and the \(3\)-class field tower \(\mathrm{F}_3^\infty(K)\) has two possible lengths \(\ell_3(K)\in\lbrace 3,4\rbrace\).
\end{theorem}

\begin{proof}
This is an immediate consequence of Table
\ref{tbl:CorrespondenceSporCC6F13},
since only the non-metabelian root \(R=B-\#4;k_1\)
with \(k_1=179\)
leads to Schur \(\sigma\)-groups \(S\)
with transfer kernel type \(\mathrm{F}.13\), \(\varkappa(K)=(3143)\),
and soluble length \(\mathrm{sl}(S)=4\),
and these have second abelian type invariants in Formula
\eqref{eqn:F13Category3Cc6b}.
\end{proof}

%--------------------------------------------------------------------------------

\begin{example}
\label{exm:F13Category3Cc6}
According to Table
\ref{tbl:ATI2cc6Spor},
the quadratic field \(K=\mathbb{Q}(\sqrt{d})\)
with discriminant \(d=-6\,863\,219\)
satisfies Formula
\eqref{eqn:F13Category3Cc6a}
of Theorem
\ref{thm:F13Category3Cc6}.
Thus, its \(3\)-class field tower has exact length \(\ell_3(K)=3\)
and relative degree \(\lbrack\mathrm{F}_3^\infty(K):K\rbrack=3^{26}\).
\end{example}

%--------------------------------------------------------------------------------

\renewcommand{\arraystretch}{1.0}

\begin{table}[ht]
\caption{Correspondence for sporadic \(M=S/S^{\prime\prime}\), \(\mathrm{cc}(M)=6\), of type \(\mathrm{F}.12\)}
\label{tbl:CorrespondenceSporCC6F12}
\begin{center}
\begin{tabular}{|c|c||c|c|c|c|c||c|c|c|c|}
\hline
   Type &  &  &  &  &  &  & \multicolumn{4}{|c|}{\(\alpha^{(2)}(S)\)} \\
 lo\((M)\)& \(j\) & \(k_1\) & \(k_3\) & lo\((S)\) & sl\((S)\) & \(\#\) & \(A_1\) & \(A_2\) & \(A_3\) & \(A_4\) \\
\hline
 \(\mathrm{F}.12\) &  &  &  &  &  &  & \multicolumn{4}{|c|}{\(\varkappa(S)=(1343)\)} \\
\hline
 \(13\) & \(47\) & \(148\) & \(22\) & \(29\) & \(3\) & \(1/45\) & \((4211)^3\) & \((4211)^3\) & \((221)^3,(211)^9\) & \((221)^3,(211)^9\) \\
   \(\) &   \(\) &    \(\) & \(39\) & \(32\) & \(4\) & \(1/45\) & \((4211)^3\) & \((4211)^3\) & \((221)^3,(211)^9\) & \((221)^3,(211)^9\) \\
   \(\) &   \(\) & \(179\) & \(22\) & \(26\) & \(3\) & \(1/45\) & \((5211)^3\) & \((4211)^3\) & \((221)^3,(211)^9\) & \((221)^3,(211)^9\) \\
   \(\) &   \(\) &    \(\) & \(39\) & \(26\) & \(3\) & \(1/45\) & \((5211)^3\) & \((4211)^3\) & \((221)^3,(211)^9\) & \((221)^3,(211)^9\) \\
 \(13\) & \(50\) & \(148\) &  \(4\) & \(32\) & \(4\) & \(1/45\) & \((4211)^3\) & \((4211)^3\) & \((221)^3,(211)^9\) & \((221)^3,(211)^9\) \\
   \(\) &   \(\) &    \(\) & \(28\) & \(35\) & \(4\) & \(1/45\) & \((4211)^3\) & \((4211)^3\) & \((221)^3,(211)^9\) & \((221)^3,(211)^9\) \\
   \(\) &   \(\) & \(179\) &  \(4\) & \(26\) & \(3\) & \(1/45\) & \((5211)^3\) & \((4211)^3\) & \((221)^3,(211)^9\) & \((221)^3,(211)^9\) \\
   \(\) &   \(\) &    \(\) & \(28\) & \(26\) & \(3\) & \(1/45\) & \((5211)^3\) & \((4211)^3\) & \((221)^3,(211)^9\) & \((221)^3,(211)^9\) \\
 \(13\) & \(55\) & \(148\) &  \(2\) & \(29\) & \(3\) & \(1/45\) & \((4211)^3\) & \((4211)^3\) & \((221)^3,(211)^9\) & \((221)^3,(211)^9\) \\
   \(\) &   \(\) &    \(\) & \(14\) & \(29\) & \(3\) & \(1/45\) & \((4211)^3\) & \((4211)^3\) & \((221)^3,(211)^9\) & \((221)^3,(211)^9\) \\
   \(\) &   \(\) & \(179\) &  \(2\) & \(26\) & \(3\) & \(1/45\) & \((5211)^3\) & \((4211)^3\) & \((221)^3,(211)^9\) & \((221)^3,(211)^9\) \\
   \(\) &   \(\) &    \(\) & \(40\) & \(26\) & \(3\) & \(1/45\) & \((5211)^3\) & \((4211)^3\) & \((221)^3,(211)^9\) & \((221)^3,(211)^9\) \\
 \(13\) & \(57\) & \(148\) & \(16\) & \(29\) & \(3\) & \(1/45\) & \((4211)^3\) & \((4211)^3\) & \((221)^3,(211)^9\) & \((221)^3,(211)^9\) \\
   \(\) &   \(\) &    \(\) & \(19\) & \(32\) & \(4\) & \(1/45\) & \((4211)^3\) & \((4211)^3\) & \((221)^3,(211)^9\) & \((221)^3,(211)^9\) \\
   \(\) &   \(\) & \(179\) & \(16\) & \(26\) & \(3\) & \(1/45\) & \((5211)^3\) & \((4211)^3\) & \((221)^3,(211)^9\) & \((221)^3,(211)^9\) \\
   \(\) &   \(\) &    \(\) & \(44\) & \(26\) & \(3\) & \(1/45\) & \((5211)^3\) & \((4211)^3\) & \((221)^3,(211)^9\) & \((221)^3,(211)^9\) \\
\hline
\end{tabular}
\end{center}
\end{table}

%--------------------------------------------------------------------------------

\begin{theorem}
\label{thm:F12Category3Cc6}
Let \(K=\mathbb{Q}(\sqrt{d})\) be an imaginary quadratic field
with fundamental discriminant \(d<0\),
elementary bicyclic \(3\)-class group \(\mathrm{Cl}_3(K)\simeq C_3\times C_3\),
capitulation type \(\mathrm{F}.12\), \(\varkappa(K)\sim (1343)\),
and first abelian type invariants
\(\alpha^{(1)}(K)=\lbrack 11;43,43,111,111\rbrack\).
If the second abelian type invariants are of third category,
they must either have the shape
\begin{equation}
\label{eqn:F12Category3Cc6a}
\alpha^{(2)}(K)=\lbrack 11;(43;3332,(5211)^3),(43;3332,(4211)^3),(111;3332,(221)^3,(211)^9)^2\rbrack
\end{equation}
and the \(3\)-class field tower \(\mathrm{F}_3^\infty(K)\) has precise length \(\ell_3(K)=3\) and relative degree \(3^{26}\) over \(K\),
or the shape
\begin{equation}
\label{eqn:F12Category3Cc6b}
\alpha^{(2)}(K)=\lbrack 11;(43;3332,(4211)^3),(43;3332,(4211)^3),(111;3332,(221)^3,(211)^9)^2\rbrack
\end{equation}
and the \(3\)-class field tower \(\mathrm{F}_3^\infty(K)\) has two possible lengths \(\ell_3(K)\in\lbrace 3,4\rbrace\).
\end{theorem}

\begin{proof}
This is an immediate consequence of Table
\ref{tbl:CorrespondenceSporCC6F12},
since only the non-metabelian root \(R=B-\#4;k_1\)
with \(k_1=148\)
leads to Schur \(\sigma\)-groups \(S\)
with transfer kernel type \(\mathrm{F}.12\), \(\varkappa(K)=(1343)\),
and soluble length \(\mathrm{sl}(S)=4\),
and these have second abelian type invariants in Formula
\eqref{eqn:F12Category3Cc6b}.
\end{proof}

%--------------------------------------------------------------------------------

\begin{example}
\label{exm:F12Category3Cc6}
According to Table
\ref{tbl:ATI2cc6Spor},
the quadratic field \(K=\mathbb{Q}(\sqrt{d})\)
with discriminant \(d=-8\,751\,215\)
satisfies Formula
\eqref{eqn:F12Category3Cc6a}
of Theorem
\ref{thm:F12Category3Cc6}.
Thus, its \(3\)-class field tower has exact length \(\ell_3(K)=3\)
and relative degree \(\lbrack\mathrm{F}_3^\infty(K):K\rbrack=3^{26}\).
However,
the quadratic field \(K=\mathbb{Q}(\sqrt{d})\)
with discriminant \(d=-423\,640\)
satisfies Formula
\eqref{eqn:F12Category3Cc6b}
of Theorem
\ref{thm:F12Category3Cc6},
and no sharp decision about the length of the \(3\)-class field tower \(\mathrm{F}_3^\infty(K)\)
can be drawn, since there are two possibilities \(\ell_3(K)\in\lbrace 3,4\rbrace\).
\end{example}

%--------------------------------------------------------------------------------

\renewcommand{\arraystretch}{1.0}

\begin{table}[ht]
\caption{Correspondence for sporadic \(M=S/S^{\prime\prime}\), \(\mathrm{cc}(M)=6\), of type \(\mathrm{F}.11\)}
\label{tbl:CorrespondenceSporCC6F11}
\begin{center}
\begin{tabular}{|c|c||c|c|c|c|c||c|c|c|c|}
\hline
   Type &  &  &  &  &  &  & \multicolumn{4}{|c|}{\(\alpha^{(2)}(S)\)} \\
 lo\((M)\)& \(j\) & \(k_1\) & \(k_3\) & lo\((S)\) & sl\((S)\) & \(\#\) & \(A_1\) & \(A_2\) & \(A_3\) & \(A_4\) \\
\hline
 \(\mathrm{F}.11\) &  &  &  &  &  &  & \multicolumn{4}{|c|}{\(\varkappa(S)=(1143)\)} \\
\hline
 \(13\) & \(40\) & \(148\) &  \(1\) & \(32\) & \(4\) & \(1/45\) & \((5211)^3\) & \((4211)^3\) & \((221)^3,(211)^9\) & \((221)^3,(211)^9\) \\
   \(\) &   \(\) &    \(\) & \(33\) & \(29\) & \(3\) & \(1/45\) & \((5211)^3\) & \((4211)^3\) & \((221)^3,(211)^9\) & \((221)^3,(211)^9\) \\
   \(\) &   \(\) & \(179\) &  \(1\) & \(29\) & \(3\) & \(1/45\) & \((5211)^3\) & \((4211)^3\) & \((221)^3,(211)^9\) & \((221)^3,(211)^9\) \\
   \(\) &   \(\) &    \(\) & \(33\) & \(29\) & \(3\) & \(1/45\) & \((5211)^3\) & \((4211)^3\) & \((221)^3,(211)^9\) & \((221)^3,(211)^9\) \\
 \(13\) & \(42\) & \(148\) & \(18\) & \(29\) & \(3\) & \(1/45\) & \((5211)^3\) & \((4211)^3\) & \((221)^3,(211)^9\) & \((221)^3,(211)^9\) \\
   \(\) &   \(\) &    \(\) & \(35\) & \(32\) & \(4\) & \(1/45\) & \((5211)^3\) & \((4211)^3\) & \((221)^3,(211)^9\) & \((221)^3,(211)^9\) \\
   \(\) &   \(\) & \(179\) & \(18\) & \(29\) & \(3\) & \(1/45\) & \((5211)^3\) & \((4211)^3\) & \((221)^3,(211)^9\) & \((221)^3,(211)^9\) \\
   \(\) &   \(\) &    \(\) & \(35\) & \(32\) & \(4\) & \(1/45\) & \((5211)^3\) & \((4211)^3\) & \((221)^3,(211)^9\) & \((221)^3,(211)^9\) \\
\hline
\end{tabular}
\end{center}
\end{table}

%--------------------------------------------------------------------------------

\begin{theorem}
\label{thm:F11Category3Cc6}
Let \(K=\mathbb{Q}(\sqrt{d})\) be an imaginary quadratic field
with fundamental discriminant \(d<0\),
elementary bicyclic \(3\)-class group \(\mathrm{Cl}_3(K)\simeq C_3\times C_3\),
capitulation type \(\mathrm{F}.11\), \(\varkappa(K)\sim (1143)\),
and first abelian type invariants
\(\alpha^{(1)}(K)=\lbrack 11;43,43,111,111\rbrack\).
If the second abelian type invariants are of third category,
they must have the shape
\begin{equation}
\label{eqn:F11Category3Cc6}
\alpha^{(2)}(K)=\lbrack 11;(43;3332,(5211)^3),(43;3332,(4211)^3),(111;3332,(221)^3,(211)^9)^2\rbrack
\end{equation}
and the \(3\)-class field tower \(\mathrm{F}_3^\infty(K)\) has two possible lengths \(\ell_3(K)\in\lbrace 3,4\rbrace\)..
\end{theorem}

\begin{proof}
This is an immediate consequence of Table
\ref{tbl:CorrespondenceSporCC6F11},
since both non-metabelian roots \(R=B-\#4;k_1\)
with \(k_1\in\lbrace 148,179\rbrace\)
lead to Schur \(\sigma\)-groups \(S\)
with transfer kernel type \(\mathrm{F}.11\), \(\varkappa(K)=(1143)\),
whose soluble length is \(\mathrm{sl}(S)=4\),
and which have the second abelian type invariants in Formula
\eqref{eqn:F11Category3Cc6}.
\end{proof}

%--------------------------------------------------------------------------------

\begin{example}
\label{exm:F11Category3Cc6}
According to Table
\ref{tbl:ATI2cc6Spor},
the quadratic fields \(K=\mathbb{Q}(\sqrt{d})\)
with discriminants \(d\in\lbrace -4\,838\,891,-8\,439\,815\rbrace\)
satisfy the assumptions of Theorem
\ref{thm:F11Category3Cc6}.
Thus, no sharp decision about the length of the \(3\)-class field tower \(\mathrm{F}_3^\infty(K)\)
can be drawn, since there are two possibilities \(\ell_3(K)\in\lbrace 3,4\rbrace\).
\end{example}

%--------------------------------------------------------------------------------

\begin{corollary}
\label{cor:SporadicCc6Category3}
The logarithmic order
of the Schur \(\sigma\)-group \(S=\mathrm{Gal}(\mathrm{F}_3^\infty(K)/K)\)
is generally \(\mathrm{lo}(S)\in\lbrace 26,29\rbrace\) for a \(3\)-class field tower of length \(\ell_3(K)=3\),
and it depends on the capitulation type \(\varkappa(K)\),
\[
\mathrm{lo}(S)\in
\begin{cases}
\lbrace 32,35\rbrace           & \text{ for type } \mathrm{F}.12, \\
\lbrace 32\rbrace  & \text{ for type } \mathrm{F}.11 \text{ or } \mathrm{F}.13,
\end{cases}
\] 
in the case of a tower of length \(\ell_3(K)=4\).
\end{corollary}

\begin{proof}
This is a common consequence of the Tables
\ref{tbl:CorrespondenceSporCC6F7},
\ref{tbl:CorrespondenceSporCC6F13},
\ref{tbl:CorrespondenceSporCC6F12},
\ref{tbl:CorrespondenceSporCC6F11}.
\end{proof}

%\newpage
%--------------------------------------------------------------------------------

\section{Historical remarks}
\label{s:History}

\noindent
The goal of this article was the investigation of the
transfer kernel type \(\mathrm{F}\)
\cite[\S\ 2.I.F, p. 36]{SoTa1934}.
Scholz and Taussky also defined
transfer kernel type \(\mathrm{D}\)
\cite[\S\ 2.I.D, p. 35]{SoTa1934} and
transfer kernel type \(\mathrm{E}\)
\cite[\S\ 2.I.E, p. 36]{SoTa1934}.
The unique Schur \(\sigma\)-groups \(S\) of type \(\mathrm{D}\) are metabelian and have
order \(\mathrm{ord}(S)=3^{5}\), soluble length \(\mathrm{sl}(S)=2\),
coclass \(\mathrm{cc}(S)=2\) and nilpotency class \(\mathrm{cl}(S)=3\).
They were identified by Scholz and Taussky
by means of annihilator ideals and Schreier polynomials
\cite{Ma2018b}.
The smallest non-metabelian Schur \(\sigma\)-groups \(S\) of type \(\mathrm{E}\) have
order \(\mathrm{ord}(S)=3^{8}\), soluble length \(\mathrm{sl}(S)=3\),
coclass \(\mathrm{cc}(S)=3\) and nilpotency class \(\mathrm{cl}(S)=5\).
They were identified by Bush and Mayer
by means of polycyclic power commutator presentations
\cite{BuMa2015}.
Their metabelianizations \(M=S/S^{\prime\prime}\)
of order \(\mathrm{ord}(S)=3^{7}\), coclass \(\mathrm{cc}(S)=2\) and nilpotency class \(\mathrm{cl}(S)=5\)
were already known to Scholz and Taussky
in terms of annihilator ideals and Schreier polynomials,
but these authors erroneously claimed that these metabelianizations \(M\)
are \(3\)-class field tower groups of imaginary quadratic number fields,
which is impossible since \(d(M)=2<3=r(M)\).

%\newpage
%--------------------------------------------------------------------------------

%--------------------------------------------------------------------------------

\end{document}